\documentclass[12pt,twoside]{amsart}
\usepackage{amssymb}

\nonstopmode

\numberwithin{equation}{section} 

\font\tengothic=eufm10 scaled\magstep 1 \font\sevengothic=eufm7
scaled\magstep 1
\newfam\gothicfam
      \textfont\gothicfam=\tengothic
      \scriptfont\gothicfam=\sevengothic

\newtheorem{theorem}{Theorem}[section]
\newtheorem{lemma}[theorem]{Lemma}
\newtheorem{proposition}[theorem]{Proposition}
\newtheorem{corollary}[theorem]{Corollary}

\theoremstyle{definition}
\newtheorem{definition}[theorem]{Definition} 
\newtheorem{remark}[theorem]{Remark}
\newtheorem{example}[theorem]{Example}

\newtheorem{notation}[theorem]{Notation}

\newtheorem{problem}[theorem]{Problem}

\newcommand{\Hom}{\operatorname{Hom}}

\newcommand{\Ext}{\operatorname{Ext}}

\newcommand{\cD}{{\mathcal D}}
\newcommand{\cE}{{\mathcal E}}
\newcommand{\cF}{{\mathcal F}}
\newcommand{\cG}{{\mathcal G}}
\newcommand{\cO}{{\mathcal O}}
\newcommand{\cL}{{\mathcal L}}

\newcommand{\cS}{{\mathcal S}}

\newcommand{\cH}{{\mathcal H}}

\newcommand {\CC}{\mathbb{C}}
\newcommand {\ZZ}{\mathbb{Z}}

\newcommand {\PP}{\mathbb{P}}

\newcommand {\mm}{{\textbf{m}}}

\def\hookright#1{\smash{ \mathop{\hookrightarrow}
    \limits^{#1}}}

\begin{document}
\title[Geometric collections and Castelnuovo-Mumford Regularity]
{Geometric collections and Castelnuovo-Mumford Regularity}

\author[L.\ Costa, R.M.\ Mir\'o-Roig]{L.\ Costa$^*$, R.M.\
Mir\'o-Roig$^{**}$}

\address{Facultat de Matem\`atiques,
Departament d'Algebra i Geometria, Gran Via de les Corts Catalanes
585, 08007 Barcelona, SPAIN } \email{costa@ub.edu, miro@ub.edu}

\date{\today}
\thanks{$^*$ Partially supported by MTM2004-00666.}
\thanks{$^{**}$ Partially supported by MTM2004-00666.}



\begin{abstract}

The paper begins by overviewing the basic facts on geometric
exceptional collections. Then, we  derive, for any coherent sheaf
$\cF$ on a smooth projective variety  with a geometric collection,
two spectral sequences: the first one abuts to $\cF$ and the
second one to its cohomology. The main goal of the paper is to
generalize Castelnuovo-Mumford regularity for coherent sheaves on
projective spaces to coherent sheaves on smooth projective
varieties $X$ with a geometric collection $\sigma $. We define the
notion of regularity of a coherent sheaf $\cF$ on $X$ with respect
to $\sigma$. We show that the basic formal properties of the
Castelnuovo-Mumford regularity of coherent sheaves over projective
spaces continue to hold in this new setting and we show that in
case of coherent sheaves on $\PP^n$ and for a suitable geometric
collection of coherent sheaves on $\PP^n$ both notions of
regularity coincide. Finally, we carefully study the regularity of
coherent sheaves on a smooth quadric hypersurface $Q_n \subset
\PP^{n+1}$ ($n$ odd) with respect to a suitable geometric
collection and we compare it with the Castelnuovo-Mumford
regularity of their extension by zero in $\PP^{n+1}$.
\end{abstract}


\maketitle

\tableofcontents


 \section{Introduction} \label{intro}

The goal of this paper is to extend the notion of
Castelnuovo-Mumford regularity for coherent sheaves on projective
spaces to coherent sheaves on smooth projective varieties with a
geometric collection with the hope to apply it to study coherent
sheaves on smooth projective varieties. As a main tool we use
geometric collections of exceptional sheaves and helix theory.

 Exceptional bundles were first considered in
\cite{DLP} by Drezet and Le Potier, where they were used to
determine  the set of triples $(r,c_1,c_2)$ such that there exists
a semistable sheaf $\cE $ on $\PP^2$ with rank $r$ and Chern
classes $c_1$ and $c_2$; and to describe moduli spaces of stable
vector bundles on $\PP^2$. In the succeeding papers \cite{G},
\cite{G1}, \cite{GR}, \cite{Kap} and \cite{Ru}, the notion of
exceptional bundle on $\PP^2$ was extended to other manifolds $X$
and even more from the category of vector bundles on $X$ to the
bounded derived category of coherent sheaves on $X$.

Any exceptional collection $(\cE_0,\cE_1,\cdots ,\cE_m)$ gives
rise to a bi-infinite collection, $\{\cE_{i}\}_{i\in\ZZ}$, called
{\em helix} and
 defined recursively by left and
right mutations (see Definition \ref{defhelix}). Helix theory was
introduced by Drezet  and Le Potier in \cite{DLP} and by
Gorodentsev and Rudakov in \cite{GR}, in connection with the
problem of constructing exceptional bundles on $\PP^n$; and helix
theory got its further progress in the succeeding papers
\cite{Bo}, \cite{CMR3}, \cite{G1},  \cite{N}, \cite{N1}, \cite{Ka}
and \cite{Ru}. Again helix theory was first developed for vector
bundles on $\PP^2$ and generalized later to any triangulated
category where for any two objects $\cE$ and $\cF$, $\Hom^{
\bullet }(\cE, \cF)$ has a structure of finite dimensional graded
vector space over $\CC$.

\vskip 2mm  In this paper, we make  an effort  to link the
abstract and general context of helix theory and exceptional
collections to concrete examples, their applications, and the
geometrical properties that we can derive. First of all, we recall
 the notion of geometric collection (see Definition
\ref{defexcellent}) introduced by Bondal and
 Polishchuk in \cite{BP}. It is well known that the length of any full
 exceptional collection of coherent sheaves $\sigma
=(\cE_0,\cE_1,\cdots ,\cE_{m})$ on a smooth projective variety $X$
of dimension $n$ is equal to the rank of the Grothendieck group
$K_0(X)$ which turns out to be greater or equal to $n+1$. We call
geometric collection any full exceptional collection of coherent
sheaves of length $n+1$. Geometric collections have nice
properties: They are automatically full strongly exceptional
collections, their strong exceptionality is preserved under
mutations, any thread of an helix associated to a geometric
collection is a full strongly exceptional collection, etc...

We address the problem of determining smooth projective varieties
with geometric collections and we prove that $\PP^n$, any quadric
hypersurface $Q_n\subset \PP^{n+1}$ ($n$ odd) and any Fano 3-fold
$X$ with $Pic(X)\cong \ZZ$ and trivial intermediate Jacobian have
geometric collections. Given a coherent sheaf $\cF$ on a smooth
projective variety $X$ with a geometric collection, we derive two
spectral sequences: A Beilinson-Kapranov type spectral sequence
which converges to $\cF$ (Theorem \ref{beil}) and an
Eilenberg-Moore type spectral sequence which abuts to the
cohomology of $\cF$ (Theorem \ref{moore}).

The existence of geometric collections allows us to generalize the
notion of Castelnuovo-Mumford regularity for coherent sheaves on
projective spaces to coherent sheaves on smooth projective
varieties with  a geometric collection (Definition \ref{new}). We
also prove that many of the main properties of the
Castelnuovo-Mumford regularity are accomplished by this new
concept. Finally, given a smooth projective variety $X
\hookright{i} \PP^n$ with a geometric collection $\sigma$ and a
coherent $\cO_X$-module $\cF$ it would be very interesting to
compare the regularity of $\cF$ with respect to $\sigma$,
$Reg_{\sigma}(\cF)$, to the Castelnuovo-Mumford regularity
$Reg^{CM}(i_*\cF)$ of its extension by zero $i_*\cF$. In the last
part of this work, we will address this problem and we will show
by carefully analyzing the cases of coherent sheaves on quadric
hypersurfaces $Q_n \subset \PP^{n+1}$ ($n$ odd) that, in general,
they are very different and we write down a formula which relates
both notions of regularity (Theorem \ref{mainquadrica}).

\vspace{3mm}

Next we outline the structure of this paper. In section 2, we
overview the basic facts on exceptional collections, geometric
collections, mutations and helix theory; we give examples to
illustrate all these concepts, and we describe the spectral
sequences that we use in the sequel to develop the theory of
regularity with respect to a geometric collection. In section 3,
we give, using a Beilinson-Kapranov type spectral sequence, the
promised definition of regularity with respect to a geometric
collection, we prove that the Castelnuovo-Mumford regularity of a
coherent sheaf $\cF$ on $\PP^n$ coincides with the regularity of
$\cF$ with respect to a suitable geometric collection of coherent
sheaves on $\PP^n$ and we show that the main basic properties of
the Castelnuovo-Mumford regularity carry over to the new setting.
In section 4, we consider odd dimensional quadric hypersurfaces
$Q_n \subset \PP^{n+1}$ and their geometric collection
$\sigma=(\cO_{Q_n}, \cO_{Q_n}(1), \cdots, \cO_{Q_n}(n-1), \Sigma
(n-1))$, where $\Sigma$ is the Spinor bundle. We compute the right
dual basis of any thread of the strict helix $\cH _{\sigma }=
\{\cE_i \}_{i \in \ZZ}$ associated to $\sigma$, we illustrate our
results on the regularity of a sheaf with respect to $\sigma $ for
the case of sheaves $\cF$ on $Q_n$ and we compare it with the
regularity of its extension by zero in the embedding $Q_n \subset
\PP^{n+1}$. This last results  show that, in general, the
regularity with respect to $\sigma $ of a coherent sheaf $\cF$ on
a smooth projective variety $X \hookright{i} \PP^m$ does not
square with the Castelnuovo-Mumford regularity of its extension by
zero $i_*\cF$ in $\PP^m$. We end the paper in \S 5 with some final
comments and questions which naturally arise from this paper.
\vspace{3mm}

\section{Geometric collections and spectral sequences}

Let $X$  be a smooth projective variety defined over the complex
numbers $\CC$ and let  $\cD=D^b({\cO}_X$-$mod)$ be the derived
category of bounded complexes of coherent sheaves of
${\cO}_X$-modules. For any pair of objects $A,B \in Obj(\cD)$ we
introduce the following notation:
\[\Hom^{\bullet}(A,B):= \bigoplus_{k \in \ZZ} \Ext^k_{\cD}(A,B) \]
and if $V^{\bullet}$ is a graded vector space and $A$ an object of
$\cD$, then the tensor product can be constructed as
\begin{equation} \label{tensor} V^{\bullet} \otimes A= \bigoplus _{\alpha} V^{\alpha}
\otimes A[-\alpha]. \end{equation} We will use the dualization
defined for graded vector spaces by the rule
\begin{equation} (V^{\times \bullet})_{p}=(V^{\bullet})_{-p}. \end{equation}
 A covariant cohomological functor $\mbox{Cov}^{\bullet}$ is called linear if for any $\beta$
satisfies
\begin{equation} \label{functor} \mbox{Cov}^{\beta}(V^{\bullet}
\otimes A)= \bigoplus _{\alpha} V^{\alpha} \otimes
\mbox{Cov}^{\beta}(A[-\alpha]) =   \bigoplus _{\alpha} V^{\alpha}
\otimes \mbox{Cov}^{\beta-\alpha}(A).\end{equation}

\begin{definition}\label{exceptcoll}
Let $X$ be a smooth projective variety.

(i) An object  $\cF \in \cD$ is {\bf exceptional} if
$\Hom^{\bullet} (\cF,\cF)$ is a 1-dimensional algebra generated by
the identity.

(ii) An ordered collection $(\cF_0,\cF_1,\cdots ,\cF_m)$  of
objects of $\cD$ is an {\bf exceptional collection} if each object
$\cF_{i}$ is exceptional and
$\Ext^{\bullet}_{\cD}(\cF_{k},\cF_{j})=0$ for $j<k$.

(iii) An exceptional collection $(\cF_0,\cF_1,\cdots ,\cF_m)$ of
objects of $\cD$ is a {\bf strongly exceptional collection} if in
addition $\Ext^{i}_{\cD}(\cF_j,\cF_k)=0$ for $i \neq 0$ and  $j
\leq k$.

(iv) \label{full}   An ordered collection of objects of $\cD$,
$(\cF_0,\cF_1,\cdots ,\cF_m)$,  is a {\bf full (strongly)
exceptional collection} if it is a (strongly) exceptional
collection $(\cF_0,\cF_1,\cdots ,\cF_m)$ and $\cF_0$, $\cF_1$,
$\cdots $ , $\cF_m$ generate the bounded derived category $\cD$.
\end{definition}

\begin{example} \label{prihirse} (1) ($\cO_{\PP^r}$, $\cO_{\PP^r} (1) $, $\cO_{\PP^r} (2)
$, $\cdots $, $\cO_{\PP^r} (r) $)  and ($\cO_{\PP^r}$,
$\Omega^1_{\PP^r} (1) $, $\Omega^2_{\PP^r} (2) $, $\cdots $,
$\Omega^r_{\PP^r} (r) $) are full strongly exceptional collections
of coherent sheaves on $\PP^r$.

(2) \label{gras} Let  $Gr(k,n)$ be the Grassmannian of
$k$-dimensional subspaces of the $n$-dimensional vector space and
let $\cS$ be the tautological $k$-dimensional bundle on $X$.
Denote by
 $\Sigma^{\alpha}\cS$  the space of the
irreducible representations of the group $GL(\cS)$ with highest
weight $\alpha=(\alpha_1, \ldots, \alpha_s)$ and
$|\alpha|=\sum_{i=1}^{s} \alpha_i$. Denote by $A(k,n)$ the set of
locally free sheaves $\Sigma^{\alpha}\cS$ on $Gr(k,n)$ where
$\alpha$ runs over Young diagrams fitting inside a $k \times
(n-k)$ rectangle. Set $\rho(k,n):= \sharp A(k,n)$. By \cite{Kap};
Proposition 2.2 (a) and Proposition 1.4, $A(k,n)$ can be totally
ordered in such a way that we obtain a full strongly exceptional
collection ($ E_1, \ldots, E_{\rho(k,n)}$) of locally free sheaves
on $Gr(k,n)$.

(3)  \label{quadrica} Let $Q_n\subset \PP^{n+1}$, $n\ge 2$, be the
quadric hypersurface. By \cite{Kap2}; Proposition 4.9, if $n$ is
even and $\Sigma _1$, $\Sigma _2$ are the Spinor bundles on $Q_n$,
then
\[(\Sigma _1(-n), \Sigma _2(-n), \cO_{Q_n}(-n+1), \cdots, \cO_{Q_n}(-1), \cO_{Q_n}) \]
is a full strongly exceptional collection of locally free sheaves
on $Q_n$; and if $n$ is odd and $\Sigma$ is the Spinor bundle on
$Q_n$, then
\[( \Sigma(-n), \cO_{Q_n}(-n+1), \cdots, \cO_{Q_n}(-1), \cO_{Q_n}) \]
is a full strongly exceptional collection of locally free sheaves
on $Q_n$.
\end{example}

\begin{remark} The existence of a full strongly exceptional collection
$(\cF_0,\cF_1,\cdots,\cF_m)$ of coherent sheaves on a smooth
projective variety $X$ imposes a rather  strong restriction on
$X$, namely that the Grothendieck group $K_0(X)=K_0(\cO _X-mod)$
is isomorphic to $\ZZ^{m+1}$.
\end{remark}

\begin{definition}
Let $X$ be a smooth projective variety and let $(A,B)$ be an
exceptional pair of objects of $\cD$. We define the {\bf left
mutation} of $B$, $L_AB$, and the {\bf right mutation} of $A$,
$R_BA$, with the aid of the following distinguished triangles in
the category $\cD$:
\begin{equation}  \label{t1} L_AB \rightarrow \Hom^{\bullet}(A,B) \otimes A \rightarrow B
\rightarrow L_AB[1]
\end{equation}
\begin{equation}  \label{t2} R_BA[-1] \rightarrow A \rightarrow \Hom^{\times \bullet}(A,B)
 \otimes B \rightarrow R_BA.
\end{equation}
\end{definition}

\begin{remark} \label{ortogonalitat}
If we apply $\Hom^{\bullet}(A, \cdot )$ to the triangle (\ref{t1})
and we apply $\Hom^{\bullet}(\cdot, B)$ to the triangle (\ref{t2})
we get the following orthogonality relations:
\[ \Hom^{\bullet}(A, L_AB)=0  \mbox{ and } \Hom^{\bullet}(R_BA, B)=0. \]
\end{remark}

\begin{definition}
Let $X$ be a smooth projective variety and let $\sigma=(\cE_0,
\cdots, \cE_n )$ be  an exceptional collection of objects of
$\cD$. A {\bf left mutation} (resp. {\bf right mutation}) of
$\sigma$ is defined as follows:  for any $1 \leq i \leq n$ a left
mutation $L_i$ replaces the $i$-th pair of consequent elements
$(\cE_{i-1}, \cE_i)$ by  $(L_{\cE_{i-1}}\cE_i, \cE_{i-1} )$ and a
right mutation $R_i$ replaces the same pair of consequent elements
$(\cE_{i-1}, \cE_i)$ by  $(\cE_i, R_{\cE_{i}}\cE_{i-1})$:
\[ L_i \sigma=L_{\cE_{i-1}} \sigma= ( \cE_0, \cdots ,L_{\cE_{i-1}}\cE_i,
\cE_{i-1}, \cdots \cE_n )\]
\[ R_i \sigma=R_{\cE_{i-1}} \sigma= ( \cE_0, \cdots ,\cE_i, R_{\cE_{i}}\cE_{i-1},
\cdots, \cE_n ).\]
\end{definition}

\begin{notation} It is convenient to agree that
 \[  R^{(j)}\cE_{i}=  R^{(j-1)}R\cE_{i}= R_{\cE_{i+j}} \cdots \cdots
  R_{\cE_{i+2}} R_{\cE_{i+1}} \cE_{i} \]
 and similar notation for compositions of left mutations. According to
 these  notations, mutations satisfy the following relations:
\[ \begin{array}{l}
L_iR_i=R_iL_i=Id \\ L_iL_j=L_jL_i \quad \mbox{for } |i-j|>1 \\
L_{i+1}L_iL_{i+1}=L_iL_{i+1}L_i \quad \mbox{for } 1 < i < n.
\end{array} \]
\end{notation}

\begin{proposition} \label{mutar} Let $X$ be a smooth projective variety
and let $\sigma=(\cE_0, \cdots, \cE_n )$ be  an exceptional
collection of objects of $\cD$. Then any mutation of $\sigma$ is
an exceptional collection and if $\sigma$ generates the category
$\cD$, then the mutated collection also generates $\cD$.
\end{proposition}
\begin{proof}
See \cite{Bo}; Assertion 2.1 and Lemma 2.2.
\end{proof}

\begin{remark} \label{notstrong} In general a mutation of a strongly exceptional
collection is not a strongly exceptional collection. For instance,
let $X= \PP^1 \times \PP^1$ be a smooth quadric surface in $\PP^3$
and denote by $\cO_X(a,b)= \cO_{\PP^1}(a) \boxtimes
\cO_{\PP^1}(b)$. By \cite{CMR}; Proposition 4.16, (see also
\cite{CMR2}
\[\sigma=(\cO_X, \cO_X(1,0), \cO_X(0,1), \cO_X(1,1) ) \]
is a full strongly exceptional collection of line bundles on $X$.
Using the exact sequence
\[ 0 \rightarrow \cO_X(-1,1) \rightarrow \Hom(\cO_X(0,1), \cO_X(1,1)) \otimes \cO_X(0,1)
\rightarrow  \cO_X(1,1) \rightarrow 0 \] we get that
$L_{\cO_X(0,1)} \cO_X(1,1)= \cO_X(-1,1)$. But, since
$\Ext^1(\cO_X(1,0), \cO_X(-1,1))=\CC$ the mutated exceptional
collection of line bundles on $X$
\[ L\sigma=(\cO_X, \cO_X(1,0), \cO_X(-1,1), \cO_X(0,1) )\]
is no more a strongly exceptional collection of line bundles on
$X$. We will come back to the problem of whether strongly
exceptionality is preserved under mutations.
\end{remark}

\begin{definition} \label{defdual} Let $X$ be a smooth projective variety.
 Given any full exceptional collection $\sigma=(\cE_0, \cdots, \cE_n )$ the
collection
\[ (L^{(n)}\cE_n, L^{(n-1)}\cE_{n-1}, \cdots, L^{(1)}\cE_1, \cE_0) \]
will be called {\bf left dual base of $\sigma$} and the collection
\[ (\cE_n, R^{(1)}\cE_{n-1}, \cdots, R^{(n)}\cE_0)\]
will be called {\bf right dual base of $\sigma$}.
\end{definition}

\begin{remark}
\label{dual} Given a full  exceptional collection   $\sigma=
(\cE_0, \cdots, \cE_n )$, its corresponding right and left dual
basis are uniquely determined up to a unique isomorphism and they
satisfy the following orthogonality conditions:
\[ \Hom^{\alpha}(R^{(j)}\cE_i, \cE_k)=0, \quad \Hom^{\alpha}(\cE_k, L^{(j)}\cE_i)=0\]
for all $\alpha, i, j$ and $k$ except
\[ \Hom^{k}(R^{(k)}\cE_{n-k}, \cE_{n-k})= \Hom^{n-k}(\cE_{n-k}, L^{(n-k)}\cE_{n-k})=\CC.\]
\end{remark}

The definition of helix and the first results about helices
appeared in \cite{DLP} and \cite{GR}. Let us recall its
definition.

\begin{definition} \label{defhelix}  Let $X$ be a smooth projective
variety. A {\bf helix of period $n+1$} is an infinite sequence $\{
\cE_i\}_{i \in \ZZ}$ of objects of $\cD$ such that for any $i \in
\ZZ$, $(\cE_i, \cdots, \cE_{i+n} )$ is an exceptional collection
and $\cE_{n+1+i}=R^{(n)}\cE_{i}$. \end{definition}

Any  exceptional collection of objects of $\cD$, $\sigma=(\cE_0,
\cdots, \cE_n )$, induces a unique helix by the rule
\begin{equation}   \cE_{n+i}=R^{(n)}\cE_{i-1} \quad \mbox{and
}\quad  \cE_{-i}=L^{(n)}\cE_{n-i+1}, \quad i>0. \end{equation}
 In that case we say that the helix is generated by $\sigma$ and that the
collection $\cH_{\sigma}:=\{ \cE_i\}_{i \in \ZZ}$ is the helix
associated to $\sigma$.  Each collection $\sigma_i=(\cE_i,
\cE_{i+1}, \cdots, \cE_{i+n} )$ is called a {\bf thread of the
helix} and it is clear that a helix is generated by any of its
thread.

\begin{example} \label{helixpn}
Let $\sigma=(\cO_{\PP^r}, \cO_{\PP^r} (1) , \cO_{\PP^r} (2) ,
\cdots , \cO_{\PP^r} (r) )$ be the full exceptional collection of
line bundles on $\PP^r=\PP(V)$ given in Example \ref{prihirse},
(1). The helix associated to $\sigma$ is given by $\cH_{\sigma}=
\{\cO_{\PP^r}(i) \}_{i \in \ZZ}$. Indeed, denote by $\cE_i=
\cO_{\PP^r}(i)$, $0 \leq i \leq r$. By definition, for any $i>0$,
$\cE_{r+i}=R^{(r)}\cE_{i-1}$ and $\cE_{-i}=L^{(r)}\cE_{r-i+1}$.
Using  the exterior powers of the Euler sequence
\begin{equation} \label{eulerwedge} 0\longrightarrow \wedge
^{k-1}T_{\PP^r}\longrightarrow \wedge ^{k}V \otimes
\cO_{\PP^r}(k)\longrightarrow \wedge ^{k}T_{\PP^r}\longrightarrow
0 \end{equation} we deduce that for any $k >0$
\[ R^{(k)} \cE_{0}=R_{\cE_k}R_{\cE_{k-1}} \cdots R_{\cE_1}\cE_0=\wedge^kT_{\PP^r}. \]
So, in particular we get $\cE_{r+1}=  R^{(r)} \cE_{0}=
\wedge^rT_{\PP^r}= \cO_{\PP^r}(r+1)$. More in general, repeating
the process and using once more the exact sequence
(\ref{eulerwedge}) we get
\[ \cE_{r+i} =R^{(r)} \cE_{i-1}=
\wedge^rT_{\PP^r}(i-1)= \cO_{\PP^r}(r+i).\] Analogously, we deduce
that $\cE_{-i}= \cO_{\PP^r}(-i)$ and hence we get that the helix
associated to $\sigma$ is given by $\cH_{\sigma}= \{\cO_{\PP^r}(i)
\}_{i \in \ZZ}$.  Notice that in this case, any thread
$\sigma_i=(\cO_{\PP^r}(i), \cO_{\PP^r} (i+1) , \cO_{\PP^r} (i+2) ,
\cdots , \cO_{\PP^r} (i+r) )$ of the helix is again a full
strongly exceptional collection of line bundles. Moreover, we have
that for any $i \in \ZZ$, \[ \cE_{i}=\cE_{i+r+1} \otimes
K_{\PP^{r}}, \] where $K_{\PP^{r}}= \cO_{\PP^r}(-r-1)$ is the
canonical line bundle.
\end{example}

This last observation is indeed a more general fact:

\begin{remark}
\label{period} Let $X$ be a smooth projective variety of dimension
$m$ and let $\sigma=(\cE_0, \cdots, \cE_n )$ be a full exceptional
collection of objects of $\cD$. Then the helix $\cH_{\sigma}$
associated to $\sigma$ has the following property of periodicity:
for any  $i \in \ZZ$,
\[ \cE_{i}= \cE_{i+n+1} \otimes K_X[m-n] \]
where $K_X$ is the canonical line bundle on $X$ and the number in
square brackets denotes the multiplicity of the shift of an object
to the left viewed as a graded complex in $\cD$.
\end{remark}

Bondal and  Polishchuk introduced in \cite{BP} the notion of
geometric collection as the exceptional  collection $\sigma $ of
objects of $\cD$ that generates a geometric helix, i.e., an helix
$\cH_{\sigma}= \{\cE_i \}_{i \in \ZZ}$ such that for any $k>0$ and
any $i \leq j$, $\Ext^k(\cE_i,\cE_j)=0$. Then they proved that
full geometric collections are exactly  full exceptional
collections of length equal to the dimension of the variety plus
one (see \cite{BP}; Proposition 3.3). We find it more convenient
to use this latter property as the definition:

\begin{definition}
\label{defexcellent}  Let $X$ be a smooth projective variety of
dimension $n$. We call  {\bf geometric collection} to any full
exceptional collection of coherent sheaves on $X$ of length $n+1$
and we call {\bf strict helix} to the helix generated by it.
\end{definition}

\begin{remark} (1)  Notice that since all full strongly
exceptional collections of coherent sheaves on $X$ have the same
length and it is equal to the $rank(K_0(X)) \geq n+1$, the length
of a geometric collection is the minimum possible.

(2) The existence of a geometric collection on a smooth variety
$X$ imposes a strong restriction on $X$, namely that its
Grothendieck group $K_0(X)\cong \ZZ^{dimX+1}$ and $X$ is forced to
be a Fano variety (see \cite{BP}; Theorem 3.4).
\end{remark}

Geometric collections have a nice behavior. For instance,
geometric collections are automatically strongly exceptional
collections of coherent sheaves and its strongly exceptionality is
preserved under mutations (Recall that, in general, strongly
exceptionality is not preserved under mutations as we have showed
in Remark \ref{notstrong}). More precisely we have:

\begin{proposition}
\label{excellent} Let $X$ be a smooth projective variety of
dimension $n$ and let $\sigma=(\cE_0, \cdots, \cE_n)$ be a
geometric collection of coherent sheaves on $X$. Then,

(i) Any mutation of the collection $\sigma$ consists also of
sheaves, i.e. complexes concentrated in the zero component of the
grading.

(ii) The collection $\sigma$ is a full strongly exceptional
collection of coherent sheaves.

(iii) Any mutation of $\sigma$ is a full strongly exceptional
collection of coherent sheaves.

(iv) Any thread $(\cE_{i}, \cE_{i+1}, \cdots, \cE_{n+i})$ of the
helix $\cH_{\sigma}$ associated to $\sigma$ is a full strongly
exceptional collection of coherent sheaves on $X$.
\end{proposition}
\begin{proof} It follows from \cite{Bo}; Assertion 9.2, Theorem
9.3 and Corollary 9.4.
\end{proof}

\begin{example}
\label{exemplesexcellents} According to Example \ref{prihirse},
  there exists a geometric
collection of coherent sheaves on $\PP^r$ and on  $Q_n \subset
\PP^{n+1}$, $n$  odd, and hence both varieties have strict
helices. On the other hand, there are no  geometric collections of
coherent sheaves on $Q_n$ for $n$  even and on
  $Gr(k,n)=Gr(k,V)$, $2\le k \le n-2$.
\end{example}

\vspace{3mm}

\begin{proposition}
\label{3folds} Let $X$ be any smooth Fano threefold with $Pic(X)
\cong \ZZ$ and trivial intermediate Jacobian. Then, $X$ has a
geometric collection.
\end{proposition}

\begin{proof} According to the classification of Fano threefolds
$X$ with $Pic(X) \cong \ZZ$ and trivial  intermediate Jacobian
(\cite{I}; Table 3.5), there exist four kinds of such manifolds:
The projective space $\PP^3$, a smooth quadric $Q_3 \subset
\PP^4$, the manifold $V_5 \subset \PP^6$ and the family of
manifolds $V_{22} \subset \PP^{12}$. The cases $X \cong \PP^3$ and
$X \cong Q_3$ follow from Example \ref{exemplesexcellents}, $(1)$
and $(2)$, respectively. The case $X \cong V_5$ is due to Orlov
(\cite{Or}) and the case
 $X\cong V_{22}$ is due to Kuznetsov (\cite{Kuz};Theorem 3).
\end{proof}

\begin{remark}
\label{exacta} Restricting ourselves to the strict helix area
$\cH_{\sigma}$, $\sigma= (\cE_0, \cdots, \cE_n)$ all the theory
can be pulled down from the triangulated category
$\cD=D^b({\cO}_X$-$mod)$ of bounded complexes of coherent sheaves
into the category of coherent sheaves Coh$(X)$. For instance, for
any $i<j$ the canonical distinguished triangles
\begin{equation}  \label{tt1} L_{\cE_i}\cE_j \rightarrow \Hom^{\bullet}(\cE_i,\cE_j)
 \otimes \cE_i \rightarrow \cE_j
\rightarrow L_{\cE_i}\cE_j[1]
\end{equation}
\begin{equation}  \label{tt2} R_{\cE_j}\cE_i[-1] \rightarrow \cE_i \rightarrow
\Hom^{\times \bullet}(\cE_i,\cE_j)
 \otimes \cE_j \rightarrow R_{\cE_j}\cE_i,
\end{equation}
turn to usual triples of coherent sheaves
\begin{equation}  \label{s1}  0 \rightarrow L_{\cE_i}\cE_j
\rightarrow \Hom(\cE_i,\cE_j)
 \otimes \cE_i \rightarrow \cE_j
\rightarrow 0
\end{equation}
\begin{equation}  \label{s2} 0 \rightarrow \cE_i \rightarrow
\Hom^{*}(\cE_i,\cE_j)
 \otimes \cE_j \rightarrow R_{\cE_j}\cE_i \rightarrow 0.
\end{equation}
So, without loss of generality, we could define the mutations
inside a strict helix in terms of  sheaves.
\end{remark}

\begin{theorem} \label{beil} {\bf (Beilinson-Kapranov type spectral sequence)}
Let $X$ be a smooth projective  variety of dimension $n$
 with a geometric collection $\sigma =(\cE_0, \cE_1,\cdots ,\cE _{n})$  and
 let ${\mathcal H}_{\sigma }=\{ \cE _{i}\} _{i\in \ZZ}$ be the  associated strict
 helix. Then for any thread $\sigma_i=(\cE_{i}, \cdots,
 \cE_{i+n})$ of the helix $\cH_{\sigma}$ and any coherent sheaf
 $\cF$on $X$ there is a spectral sequence with $E_1$-term
 \begin{equation}\label{sucespectral0}
^{i}E_1^{pq}=\Ext^q(R^{(-p)}\cE _{i+n+p},{\cF})\otimes \cE
_{i+p+n}\end{equation} situated in the square $0\le q\le n $, $-n
\le p \le 0$ which converges to $$ E^{i}_{\infty }=
\begin{cases} {\cF} \mbox{ for } i=0 \\ 0 \mbox{ for } i\ne 0.\end{cases}$$
\end{theorem}
\begin{proof}
First of all notice that since the helix $\cH_{\sigma}$ is strict,
the thread $\sigma_i=(\cE_{i}, \cdots, \cE_{i+n})$ also generates
the category $\cD$. We write $V^{\bullet}_k$ for the graded vector
spaces
\[ V^{\bullet}_k= \Hom^{\bullet}(R^{(n-k)}\cE_{k+i}, \cF) \]
and we consider the complex
\[L^{\bullet}: 0 \rightarrow V_0^{\bullet} \otimes \cE_{i}
 \rightarrow V_1^{\bullet} \otimes \cE_{i+1} \rightarrow \cdots
 \rightarrow V_{n-1}^{\bullet} \otimes \cE_{i+n-1} \rightarrow V_n^{\bullet} \otimes \cE_{i+n}
 \rightarrow 0   \]
 where the tensor product is defined as in (\ref{tensor}). The
 right mutations produce a canonical right Postnikov system of the
 complex $L^{\bullet}$, which naturally identifies $\cF$ with the
 canonical right convolution of this complex. Then, for an
 arbitrary linear covariant cohomological functor
 $\Phi^{\bullet}$, there exists an spectral sequence with
 $E_1$-term
 \[ ^{i}E_1^{pq}=\Phi^q(L^p)\]
 situated in the square $0 \leq p,q \leq n$ and converging to
  $\Phi^{p+q}(\cF)$ (see \cite{Kap2}; 1.5). Since
 $\Phi^{\bullet}$ is a linear functor, it follows from
 (\ref{functor}) that
\begin{equation} \label{espectralgral} \Phi^q(L^p)= \Phi^q(V_p^{\bullet} \otimes \cE_{i+p})=
\bigoplus_{l} V_p^l \otimes \Phi^{q-l}(\cE_{i+p})=
\bigoplus_{\alpha+\beta=q} V_p^{\alpha} \otimes
\Phi^{\beta}(\cE_{i+p}). \end{equation} In particular, if we
consider the covariant linear cohomology functor which takes a
complex to its cohomology sheaves and acts identically on
sheaves, i.e.
\[ \Phi^{\beta}(\cF)= \begin{cases} \cF \mbox{ for } \beta=0 \\ 0 \mbox{ for } \beta \neq 0
 \end{cases}\]
 on any  sheaf $\cF$, in the square $0 \leq p,q \leq n$, we get
 \[ ^{i}E_1^{pq}=V_p^{q} \otimes
\cE_{i+p} = \Ext^{q}(R^{(n-p)}\cE_{p+i}, \cF) \otimes  \cE_{i+p}\]
where the last equality follows from the definition of
$V_k^{\bullet}$. Finally, if we call $p'=p-n$ we get the spectral
sequence with $E_1$-term
 \[ ^{i}E_1^{p'q}=\Ext^{q}(R^{(-p')}\cE_{p'+n+i}, \cF) \otimes  \cE_{p'+n+i}\]
 in the square $0 \leq q \leq n$, $-n \leq p' \leq 0$ and which
 converges to
 $$ E^{i}_{\infty }=
\begin{cases} {\cF} \mbox{ for } i=0 \\ 0 \mbox{ for } i\ne
0.\end{cases}$$
\end{proof}

\begin{theorem} \label{moore} {\bf (Eilenberg-Moore type spectral sequence)}
Let $X$ be a smooth projective  variety of dimension $n$
 with a geometric collection $\sigma =(\cE_0, \cE_1,\cdots ,\cE _{n})$  and
 let ${\mathcal H}_{\sigma }=\{ \cE _{i}\} _{i\in \ZZ}$ be the  associated strict
 helix. Then for any thread $\sigma_i=(\cE_{i}, \cdots,
 \cE_{i+n})$ of the helix $\cH_{\sigma}$ and any pair of coherent
 sheaves $\cF$, $\cG$ on $X$ there is a spectral sequence with $E_1$-term
 \begin{equation}\label{sucespectral2}
^{i}E_1^{pq}=\bigoplus_{ \alpha+\beta=q} \Ext^{\alpha}(R^{(-p)}\cE
_{i+n+p},{\cF})\otimes \Ext^{\beta}(\cG,\cE
_{i+p+n})\end{equation} situated in the square $0\le q\le n $, $-n
\le p \le 0$ which converges to $$ E^{i}_{\infty }=
\begin{cases} {\Ext^{p+q}(\cG, \cF)} \mbox{ for } i=0 \\ 0 \mbox{ for } i\ne 0.\end{cases}$$
\end{theorem}
\begin{proof} We follow step by step  the proof of Theorem
\ref{beil} but in this case in (\ref{espectralgral}) we take as
$\Phi^{\bullet}$ the covariant functor $\Hom^{\bullet}(\cG,
\cdot)$ and in such a way we get the spectral sequence  \[
^{i}E_1^{pq}=\bigoplus_{ \alpha+\beta=q} \Ext^{\alpha}(R^{(-p)}\cE
_{i+n+p},{\cF})\otimes \Ext^{\beta}(\cG,\cE _{i+p+n}) \] situated
in the square $0\le q\le n $, $-n \le p \le 0$ which converges to
$\Phi^{p+q}(\cF)=\Ext^{p+q}(\cG, \cF)$.
\end{proof}

\vspace{3mm} We will end this section with technical results that
will be used in next sections.

\begin{lemma}
\label{tecnic} Let $X$ be a smooth projective variety of dimension
$n$ and let $\sigma=(\cE_0, \cdots, \cE_n)$ be a geometric
collection of coherent sheaves on $X$. For any $i < j$ and any
invertible sheaf $\cF$, it holds:

(a) $(L_{\cE_i}\cE_j)^*=R_{\cE_i^*} \cE_{j}^*$;

(b) $(R_{\cE_j}\cE_i)^*=L_{\cE_j^*} \cE_{i}^*$;

(c) $(R_{\cE_j}\cE_i)\otimes \cF \cong R_{\cE_j \otimes \cF}(\cE_i
\otimes \cF)$ and  $(L_{\cE_i}\cE_j)\otimes \cF \cong L_{\cE_i
\otimes \cF}(\cE_j \otimes \cF)$.
\end{lemma}
\begin{proof} (a) According to Remark \ref{exacta}; (\ref{s1}),
$L_{\cE_i}\cE_j$ is given by the aid of the exact sequence
\begin{equation}  \label{e1}  0 \rightarrow L_{\cE_i}\cE_j
\rightarrow \Hom(\cE_i,\cE_j)
 \otimes \cE_i \rightarrow \cE_j
\rightarrow 0.\end{equation} Dualizing this exact sequence we get
\[ 0 \rightarrow \cE_j^* \rightarrow \Hom(\cE_j^*,\cE_i^*)
 \otimes \cE_i^* \cong \Hom^*(\cE_i,\cE_j)
 \otimes \cE_i^* \rightarrow (L_{\cE_i}\cE_j)^*\rightarrow 0 \]
which according to Remark \ref{exacta}; (\ref{s2}) gives that
$(L_{\cE_i}\cE_j)^*=R_{\cE_i^*} \cE_{j}^*$.

(b) The proof is analogous to the proof of $(a)$.

(c) Since $\Hom(\cE_i, \cE_j) \cong \Hom(\cE_i \otimes  \cF, \cE_j
\otimes \cF )$, it is enough to tensor by $\cF$ the exact
sequences \[ 0 \rightarrow L_{\cE_i}\cE_j \rightarrow
\Hom(\cE_i,\cE_j)
 \otimes \cE_i \rightarrow \cE_j
\rightarrow 0 \] and \[ 0 \rightarrow \cE_i \rightarrow
\Hom^{*}(\cE_i,\cE_j)
 \otimes \cE_j \rightarrow R_{\cE_j}\cE_i \rightarrow 0. \]
\end{proof}

\begin{corollary}
\label{corotecnic} Let $X$ be a smooth projective variety of
dimension $n$ with canonical line bundle $K$ and let
$\sigma=(\cE_0, \cdots, \cE_n)$ be a geometric collection of
coherent sheaves on $X$. Assume that $\tau=(\cF_0, \cdots, \cF_n)$
is the right dual base of $\sigma$. Then, for any integer
$\lambda$, the right dual base of
 $\sigma_{\lambda(n+1)}=(\cE_{\lambda(n+1)},\cE_{\lambda(n+1)+1},
\cdots, \cE_{\lambda(n+1)+n})$ is \[\tau_{\lambda(n+1)}=(\cF_0
\otimes (K^*)^{\otimes \lambda}, \cF_1 \otimes (K^*)^{\otimes
\lambda}, \cdots, \cF_n \otimes (K^*)^{\otimes \lambda}). \]
\end{corollary}
\begin{proof} By definition of right dual base we have
\[ \cF_j= R^{(j)}\cE_{n-j} \quad \quad 0 \leq j \leq n.\]
On the other hand, by Remark \ref{period}
\[ \cE_{\lambda(n+1)+i}= \cE_i \otimes (K^*)^{\otimes \lambda}. \]
Therefore,  applying Lemma \ref{tecnic}, we get
\[ R^{(j)}\cE_{\lambda(n+1)+n-j}= R^{(j)}(\cE_{n-j} \otimes (K^*)^{\otimes \lambda})
= (R^{(j)}\cE_{n-j}) \otimes (K^*)^{\otimes \lambda}= \cF_j
\otimes (K^*)^{\otimes \lambda}. \]
\end{proof}

\section{Regularity with respect to geometric collections: definition and properties}

In this section, using the so-called Beilinson-Kapranov spectral
sequence, we generalize the notion of Castelnuovo-Mumford
regularity for coherent sheaves on a projective space to coherent
sheaves on a smooth projective variety with a geometric collection
of coherent sheaves. We establish for coherent sheaves on $\PP^n$
the agreement of the new definition of regularity with the old one
and we prove that many formal properties of Castelnuovo-Mumford
regularity continue to hold in our more general setup.

\vskip 2mm Let $X$ be a smooth projective variety of dimension $n$
and let $\sigma =(\cE_0, \cE_1,\cdots ,\cE _{n})$ be a geometric
collection on $X$. Associated to $\sigma $ we have a strict helix
${\mathcal H}_{\sigma }=\{ \cE_ {i} \}_{i\in \ZZ}$; and for any
collection $\sigma _i=(\cE _{i},\cE _{i+1},\cdots ,\cE _{i+n})$ of
$n+1$ subsequent sheaves (i.e. for any thread of the helix
${\mathcal H}_{\sigma }$) and any coherent ${\cO}_{X}$-module
${\cF}$ a spectral sequence (See Theorem \ref{beil})

\begin{equation}\label{sucespectral}
^{i}E_1^{pq}=\Ext^q(R^{(-p)}\cE _{i+p+n},{\cF})\otimes \cE
_{i+p+n}
\end{equation} situated in the square $0\le q\le n $, $-n
\le p \le 0$ which converges to $$ E^{i}_{\infty }=
\begin{cases} {\cF} \mbox{ for } i=0 \\ 0 \mbox{ for } i\ne 0.\end{cases}$$

\begin{definition} \label{new} Let $X$ be a smooth projective  variety of dimension $n$
 with a geometric collection $\sigma =(\cE_0, \cE_1,\cdots ,\cE _{n})$ and
 let ${\cF}$ be a coherent ${\cO}_{X}$-module. We say that
 ${\cF}$ is {\bf $m$-regular with respect to $\sigma$}
 if $\Ext^q(R^{(-p)}\cE _{-m+p},{\cF})=0$ for $q>0$ and $-n\le p \le
 0$.
\end{definition}

\vskip 2mm So, ${\cF}$ is $m$-regular with respect to $\sigma $ if
$^{-n-m}E_1^{pq}=0$ for $q>0$ in (\ref{sucespectral}). In
particular, if ${\cF}$ is $m$-regular with respect to $\sigma$ the
spectral sequence $^{-n-m}E_1^{pq}$ collapses at $E_2$  and we get
the following exact sequence:

\begin{equation} \label{re}
0 \longrightarrow {\cL}_{-n} \longrightarrow \cdots
\longrightarrow {\cL}_{-1} \longrightarrow {\cL}_0 \longrightarrow
{\cF} \longrightarrow 0
\end{equation}
where ${\cL}_p=H^0(X,(R^{(-p)}\cE_{-m+p})^{*}\otimes {\cF})\otimes
\cE_{-m+p}$ for $-n\le p \le 0$.

\vspace{3mm}

\begin{definition} Let $X$ be a smooth projective  variety of dimension $n$
 with a geometric collection $\sigma =(\cE_0, \cE_1,\cdots ,\cE _{n})$  and
 let ${\cF}$ be a coherent ${\cO}_{X}$-module. We define the
 {\bf regularity of $\cF$ with respect to $\sigma$}, $Reg_{\sigma }(\cF)$, as the least integer $m$ such
 that $\cF$ is $m$-regular with respect to $\sigma $. We set $Reg_{\sigma }(\cF)=-\infty $ if there is no such
integer.
\end{definition}

\begin{remark}
It would be nice to characterize the sheaves $\cF$ on $X$ with
$Reg_{\sigma }(\cF)=-\infty $.
\end{remark}

\begin{example} \label{ex313} Let $V$ be a $\CC $-vector space of dimension
$n+1$ and set $\PP^n=\PP(V)$. We consider the geometric collection
$\sigma =(\cO_{\PP^n}, \cO_{\PP^n}(1),\cdots ,\cO_{\PP^n}(n))$ on
$\PP^n$ and the associated strict helix ${\mathcal H}_{\sigma
}=\{\cO_{\PP^n}(i)\}_{i\in \ZZ}$ (see Example \ref{helixpn}).
Using the exterior powers of the Euler sequence
$$ 0\longrightarrow \wedge ^{k-1}T_{\PP^n}\longrightarrow \wedge
^{k}V \otimes \cO_{\PP^n}(k)\longrightarrow \wedge
^{k}T_{\PP^n}\longrightarrow 0$$   we compute the right dual basis
of any thread $\sigma _{i} =(\cO_{\PP^n}(i),
\cO_{\PP^n}(i+1),\cdots ,\cO_{\PP^n}(i+n))$ of the helix
${\mathcal H} _{\sigma }$ and we get
$$(\cO_{\PP^n}(n+i), R^{(1)}\cO_{\PP^n}(i+n-1),\cdots
,R^{(j)}\cO_{\PP^n}(i+n-j),\cdots ,R^{(n)}\cO_{\PP^n}(i))$$
$$=(\cO_{\PP^n}(n+i), T_{\PP^n}(i+n-1),\cdots
,\wedge^{j}T_{\PP^n}(i+n-j),\cdots ,\wedge^{n}T_{\PP^n}(i)).$$

Therefore, for any coherent sheaf $\cF$ on $\PP^n$ our definition
reduces to say: $\cF$ is $m$-regular with respect to $\sigma$ if
$\Ext^q(\wedge^{-p}T(-m+p),\cF)=H^{q}(\PP^n,\Omega
^{-p}(m-p)\otimes \cF)=0$ for all $q>0$ and all $p$, $-n\le p\le
0$.
\end{example}

We  now compute the regularity with respect to $\sigma =(\cE_0,
\cE_1,\cdots ,\cE _{n})$ of the sheaves $\cE _{i}$.

\begin{proposition} \label{regEi} Let $X$ be a smooth projective  variety of dimension $n$
 with a geometric collection $\sigma =(\cE_0, \cE_1,\cdots ,\cE _{n})$  and
 let ${\mathcal H}_{\sigma }=\{ \cE _{i}\} _{i\in \ZZ}$ be the
 associated strict
 helix. Then, for any $i\in \ZZ$, $Reg_{\sigma }(\cE _{i})=-i$.
\end{proposition}
\begin{proof} First of all we will see that $Reg_{\sigma }(\cE _{i})\le
-i$. By Remark \ref{ortogonalitat}, we have $$\Ext^q(R^{(-p)}\cE
_{i+p},\cE_{i})=\Ext^q(R_{\cE _{i}} \cdots R_{\cE _{i+p+2}}R_{\cE
_{i+p+1}} \cE _{i+p},\cE_{i})=0$$ for $q>0$ and $-n\le p\le 0$.
So, $\cE _{i}$ is $(-i)$-regular with respect to $\sigma $ or,
equivalently, $Reg_{\sigma }(\cE _{i})\le -i$.

Let us now see that $\cE _{i}$ is not $(-i-1)$-regular with
respect to $\sigma $. To this end, it is enough to see that there
is $ q>0$ and there is $p$, $ -n\le p \le 0$, such that
$$\Ext^q(R^{(-p)}\cE _{i+1+p},\cE_{i})\neq 0.$$ To prove it, we
write $i=\alpha n +j$ with  $0\le j < n$, $\alpha \in \ZZ$, we
consider the thread $\sigma _{\alpha n}=(\cE _{\alpha n},\cE
_{\alpha n+1},\cdots ,\cE _{\alpha n+j}=\cE _{i},\cdots ,\cE
_{\alpha n+n})$ and we construct its right dual basis (see
Definition \ref{defdual}) $$(\cE _{\alpha n+n},R^{(1)}\cE _{\alpha
n+n-1},\cdots ,R^{(n-j)}\cE _{\alpha n+j},\cdots ,R^{(n)}\cE
_{\alpha n}).$$ It follows from Remark \ref{dual} that
 $$\Ext^{n-j}(R^{(n-j)}\cE
_{\alpha n+j}, \cE _{\alpha n+j})= \Ext^{n-j}(R_{\cE _{\alpha
n+n}}R_{\cE _{\alpha n+n-1}}\cdots R_{\cE _{\alpha n+j+1}}\cE
_{\alpha n+j}, \cE _{\alpha n+j})=\CC.$$

We consider the exact sequence $$0 \longrightarrow  R_{\cE
_{\alpha n+n-1}} \cdots R_{\cE _{\alpha n+j+1}} \cE _{\alpha n+j}
\longrightarrow \Hom^*(R_{\cE _{\alpha n+n-1}} \cdots R_{\cE
_{\alpha n+j+1}} \cE _{\alpha n+j}, \cE _{\alpha n+n})\otimes \cE
_{\alpha n+n} $$ $$\longrightarrow R_{\cE _{\alpha n+n}} R_{\cE
_{\alpha n+n-1}} \cdots R_{\cE _{\alpha n+j+1}} \cE _{\alpha n+j}
\longrightarrow 0$$ and we apply the contravariant functor $\Hom
(.,\cE _{\alpha n+j})$. Since $\Ext^q(\cE _{\alpha
n+n},\cE_{\alpha n+j})= 0$ for $q>0$, we get $$\Ext^{n-j-1}(R_{\cE
_{\alpha n+n-1}} \cdots R_{\cE _{\alpha n+j+1}} \cE _{\alpha
n+j},\cE_{\alpha n+j})=\CC.$$

We repeat the process using the consequent right mutations and we
get
$$\Ext^{n-j-k}(R_{\cE _{\alpha n+n-k}} \cdots R_{\cE _{\alpha n+j+1}} \cE
_{\alpha n+j},\cE_{\alpha n+j})=\CC$$ for $0\le k\le n-1-j$. In
particular,
$$\Ext^1(R_{\cE _{i+1}}\cE _{i},\cE_{i})=\CC$$ which implies that
$\cE _{i}$ is not $(-i-1)$-regular and we conclude that
$Reg_{\sigma }(\cE _{i})=-i$.
\end{proof}

In \cite{Mum}, Lecture 14, D. Mumford defined the notion of
regularity for a coherent sheaf over a projective space. Let us
recall it.

\vspace{3mm}

\begin{definition} \label{defCM} A coherent sheaf $\cF $ on $\PP^n$ is said to
be {\bf $m$-regular in the sense of Castelnuovo-Mumford} if
$H^{i}(\PP^n,\cF (m-i))=0$ for $i>0$. We define the
Castelnuovo-Mumford regularity of $\cF$, $Reg^{CM}(\cF)$, as  the
least integer $m$ such that $\cF$ is $m$-regular.
\end{definition}

Notice that such an $m$ always exists by the ampleness of
$\cO_{\PP^n}(1)$ (\cite{H}; Chap. III, Proposition 5.3).

\vspace{3mm}

Let us now establish for coherent sheaves on $\PP^n$ the agreement
of regularity definition in the sense of Definition \ref{new} with
Castelnuovo-Mumford definition.

\begin{proposition} \label{CM}  A coherent sheaf $\cF $ on $\PP^n=\PP (V)$ is
 $m$-regular in the sense of Castelnuovo-Mumford if and only if it is
 $m$-regular with respect to the geometric collection
$\sigma =(\cO_{\PP^n}, \cO_{\PP^n}(1),\cdots ,\cO_{\PP^n}(n))$.
Hence, we have
 $$ Reg_{\sigma }(\cF)=Reg^{CM}(\cF).$$
\end{proposition}

\begin{proof} According to Definitions \ref{new} and \ref{defCM}, and Example \ref{ex313} we
have to see that
\begin{equation}\label{old1}
H^q(\PP^n,\cF(m-q))=0 \mbox{ for all } q>0
\end{equation}
if and only if
\begin{equation}\label{new1}
H^q(\PP^n,\cF\otimes \Omega _{\PP^n}^p(m+p))=0 \mbox{ for all }
q>0 \mbox{ and } 0\le p \le n.
\end{equation}

Let us first see that (\ref{old1}) implies (\ref{new1}). By
\cite{Mum}, Lecture 14, the equalities (\ref{old1}) are equivalent
to
\begin{equation}\label{old2}
H^q(\PP^n,\cF(t))=0 \mbox{ for all } q>0 \mbox{ and } t\ge m-q.
\end{equation}

Since $\Omega _{\PP^n}^n(n)\cong \cO_{\PP^n}(-1)$, we deduce from
(\ref{old2}) that
\begin{equation}\label{new2}
H^q(\PP^n,\cF\otimes \Omega _{\PP^n}^n(n+t))=0 \mbox{ for all }
q>0 \mbox{ and all } t\ge m-q+1.
\end{equation}

Using the exact cohomology sequence
\begin{equation}
\cdots \longrightarrow H^q(\PP^n,V\otimes \cF(t))\longrightarrow
H^q(\PP^n, \Omega ^{n-1}_{\PP^n}(n-1)\otimes \cF(t+1))
\longrightarrow H^{q+1}(\PP^n, \Omega ^{n}_{\PP^n}(n)\otimes
\cF(t))\longrightarrow \cdots
\end{equation}
associated to
$$0 \longrightarrow \Omega ^n_{\PP^n}(n)\longrightarrow  \wedge ^nV\otimes \cO_{\PP^n}
\longrightarrow
 \Omega ^{n-1}_{\PP^n}(n)\longrightarrow 0$$
 together with the equalities (\ref{new2}) and (\ref{old2}) we
 deduce
\begin{equation}\label{new3}
H^q(\PP^n,\cF\otimes \Omega _{\PP^n}^{n-1}(n-1+t))=0 \mbox{ for
all } q>0 \mbox{ and all } t\ge m-q+1.
\end{equation}

Going on using the exact cohomology sequence
\begin{equation}
\cdots \longrightarrow H^q(\PP^n,\wedge ^{i+1}V\otimes
\cF(t))\longrightarrow H^q(\PP^n, \Omega ^{i}_{\PP^n}(i)\otimes
\cF(t+1)) \longrightarrow H^{q+1}(\PP^n, \Omega
^{i+1}_{\PP^n}(i+1)\otimes \cF(t))\longrightarrow \cdots
\end{equation}
associated to
$$0 \longrightarrow \Omega ^{i+1}_{\PP^n}(i+1)\longrightarrow  \wedge ^{i+1}V\otimes \cO_{\PP^n}
\longrightarrow
 \Omega ^{i}_{\PP^n}(i+1)\longrightarrow 0$$
we deduce
\begin{equation}\label{new3}
H^q(\PP^n,\cF\otimes \Omega _{\PP^n}^{i}(i+t))=0 \mbox{ for all }
q>0 \mbox{ and all } t\ge m-q+1
\end{equation}
which obviously implies (\ref{new1}). Let us prove the converse.
To this end we consider the Eilenberg-Moore type spectral sequence
(see Theorem \ref{moore})

$$ E_1^{pq}=\bigoplus _{\alpha + \beta =q} \Big(\Ext^{\alpha
}(R^{(-p)}\cE _{i+n+p},\cF )\otimes \Ext^{\beta }(\cG,\cE
_{i+n+p})\Big) \quad 0\le q \le n, \quad -n\le p \le 0$$

$$d_r^{pq}:E_r^{pq} \longrightarrow E_{r}^{p+r,q-r+1}, \quad \quad
E_{\infty }\Rightarrow \Ext^{p+q}(\cG,\cF) $$ and we apply it to
the case $i=-m-n$ and $\cG=\cO _{\PP^{n}}(t-m)$. So, we have

$$ E_1^{pq}=\bigoplus _{\alpha + \beta =q} \Big(H^{\alpha
}(\PP^n,\Omega_{\PP^n}^{(-p)}\otimes \cF )\otimes H^{\beta
}(\PP^n,\cO_{\PP^n}(p-t))\Big) \quad 0\le q \le n, \quad -n\le p
\le 0.$$

By (\ref{new1}), $H^{\alpha }(\PP^n,\Omega_{\PP^n}^{(-p)}\otimes
\cF )=0$ for all $\alpha >0$ and all $p$, $-n\le p \le 0$ and by
Bott's formulas a non-zero cohomology group of line bundles on
$\PP^n$, $H^{\beta }(\PP^n,\cO_{\PP^n}(p-t))$, corresponds only to
$\beta =0$ or $n$ and in the latter case $p-t\le -n-1$ (i.e. $p\le
t-n-1$). Thus, $E_1^{pq}\ne 0$ forces $p+q\le t-1$. Therefore,
$E_{\infty}^{pq}=0$ for $p+q\ge t$ and so
$0=\Ext^t(\cO_{\PP^n}(t-m),\cF)=H^{t}(\PP^n,\cF(m-t))$ for $t\ge
1$ or, equivalently, $\cF$ is $m$-regular in the sense of
Castelnuovo-Mumford.
\end{proof}

\vspace{3mm}

Let us now prove that the main formal properties of
Castelnuovo-Mumford regularity over projective spaces remain to be
true in the new setting.

\vspace{3mm}

\begin{proposition}\label{sup}  Let $X$ be a smooth projective  variety of dimension $n$
 with a geometric collection $\sigma =(\cE_0, \cE_1,\cdots ,\cE _{n})$  and
 let ${\cF}$ be a coherent ${\cO}_{X}$-module. If $\cF $ is
$m$-regular with respect to $\sigma $ then  the canonical map
$\Hom(\cE_{-m}, \cF) \otimes \cE_{-m} \twoheadrightarrow \cF$ is
surjective and $\cF$ is $k$-regular with respect to $\sigma $ for
any $k\ge m$ as well.
\end{proposition}
\begin{proof} The first assertion immediately follows from the
exact sequence
 (\ref{re}) of $\cF$. To prove the second assertion,  it is enough to check it for
$k=m+1$. Since $\cF$ is $m$-regular with respect to $\sigma $ we
have
\begin{equation}\label{mregul} \Ext^q(R^{(-p)}\cE_{-m+p},\cF)=0
\mbox{ for all } q>0 \mbox{ and all } -n\le p\le 0.
\end{equation}

In order to prove that $\cF  $ is $(m+1)$-regular with respect to
$\sigma $
 we have to check that
\begin{equation}\label{m1regul} \Ext^q(R^{(-p)}\cE_{-m-1+p},\cF)=0
\mbox{ for all } q>0 \mbox{ and all } -n\le p\le 0.
\end{equation}

To this end we apply the contravariant functor $\Hom(.,\cF)$ to
the exact sequence

$$ 0\longrightarrow \cE _{-m-1} \longrightarrow \Hom^*(\cE
_{-m-1},\cE _{-m})\otimes \cE _{-m}\longrightarrow R_{\cE
_{-m}}\cE _{-m-1}\longrightarrow  0$$ and we get the exact
sequence

$$\cdots \longrightarrow \Ext ^q(R_{\cE _{-m}}\cE
_{-m-1},\cF)\longrightarrow \Hom^*(\cE _{-m-1},\cE _{-m})\otimes
\Ext ^q(\cE _{-m},\cF)\longrightarrow $$ $$\Ext ^q(\cE
_{-m-1},\cF)\longrightarrow \Ext ^{q+1}(R_{\cE _{-m}}\cE
_{-m-1},\cF)\longrightarrow \cdots .$$

Since by (\ref{mregul}), $\Ext ^q(\cE _{-m},\cF)=\Ext ^q(R_{\cE
_{-m}}\cE _{-m-1},\cF)=0$ for all $q>0$, we get that $\Ext ^q(\cE
_{-m-1},\cF)=0$ for all $q>0$. Using again (\ref{mregul}) and the
exact sequence

$$ 0\longrightarrow R_{\cE _{-m-1}}\cE _{-m-2} \longrightarrow
\Hom^*(R_{\cE _{-m-1}}\cE _{-m-2},\cE _{-m})\otimes \cE
_{-m}\longrightarrow R_{\cE _{-m}}R_{\cE _{-m-1}}\cE
_{-m-2}\longrightarrow 0$$ we get that $\Ext ^q(R_{\cE
_{-m-1}}\cE_{-m-2},\cF)=0$ for all $q>0$.

Going on and using the consequent right mutations, we get for all
$i$, $1\le i \le n-1$,$$\Ext ^q(R_{\cE _{-m-1}}R_{\cE
_{-m-2}}\cdots R_{\cE _{-m-i}}\cE_{-m-i-1},\cF)=0 \mbox{ for all }
q>0.$$ Hence, it only remains to see that $$\Ext ^q(R_{\cE
_{-m-1}}R_{\cE _{-m-2}}\cdots R_{\cE _{-m-n}}\cE_{-m-n-1},\cF)=0
\mbox{ for all } q>0.$$

The vanishing of these last $\Ext$'s groups follows again from
(\ref{mregul}) taking into account that, by Definition
\ref{defdual} \[ R_{\cE _{-m-1}}R_{\cE _{-m-2}}\cdots
 R_{\cE_{-m-n}}\cE_{-m-n-1}= R^{(n)} \cE_{-m-n-1}=\cE _{-m}. \]
\end{proof}

\vspace{3mm}

\begin{proposition}
\label{propietats} Let $X$ be a smooth projective  variety of
dimension $n$
 with a geometric collection $\sigma =(\cE_0, \cE_1,\cdots ,\cE _{n})$, let $\cF$
 and $\cG$ be coherent ${\cO}_{X}$-modules and
 let
\begin{equation} \label{regsuc} 0  \rightarrow \cF_1 \rightarrow \cF_2 \rightarrow \cF_3
\rightarrow 0\end{equation} be an exact sequence of coherent
${\cO}_{X}$-modules. Then,

 \begin{itemize}
 \item[(a)] $Reg_{\sigma }(\cF_2) \leq max\{Reg_{\sigma }(\cF_1),Reg_{\sigma
 }(\cF_3)\}$.
  \item[(b)] $Reg_{\sigma }(\cF \oplus \cG)=max\{Reg_{\sigma }(\cF),Reg_{\sigma
  }(\cG)\}$.
 \end{itemize}
\end{proposition}
\begin{proof}
(a) Let $m= max\{Reg_{\sigma }(\cF_1),Reg_{\sigma
 }(\cF_3)\}$. Since, by Proposition \ref{sup}, $\cF_1 $ and $\cF_3 $ are
both $m$-regular with respect to $\sigma $ considering the long
exact sequence
\[ \cdots \rightarrow \Ext^q(R^{(-p)}\cE _{-m+p},{\cF_1}) \rightarrow
 \Ext^q(R^{(-p)}\cE _{-m+p},{\cF_2}) \rightarrow \Ext^q(R^{(-p)}\cE _{-m+p},{\cF_3})
 \rightarrow \cdots \]
associated to (\ref{regsuc}) we get
$\Ext^q(R^{(-p)}\cE_{-m+p},{\cF_2})=0$ for any $q>0$ and $-n \leq
p \leq 0$, which implies that $Reg_{\sigma}(\cF_2) \leq m$.

(b) It is enough to see that if $\cF$ is $m$-regular with respect
to $\sigma $ and $\cG$ is $s$-regular with respect to $\sigma $
then $\cF \oplus \cG$ is $t=max(s,m)$-regular with respect to $
\sigma $. Since, by Proposition \ref{sup}, $\cF $ and $\cG $ are
both $t$-regular with respect to $\sigma $ and, moreover,  the
functor $\Ext $'s is  additive we get $$\Ext^q(R^{(-p)}\cE
_{-t+p},{\cF}\oplus {G })=\Ext^q(R^{(-p)}\cE
_{-t+p},{\cF})=\Ext^q(R^{(-p)}\cE _{-t+p},{\cG })=0$$ for $q>0$
and $-n\le p \le
 0$. Therefore, $\cF \oplus \cG$ is $t$-regular with
respect to $ \sigma $.
\end{proof}

\vspace{3mm}

The following Example will show that in Proposition
\ref{propietats}, $(a)$ we can have strict inequality.

\vspace{3mm}

\begin{example}
We consider the geometric collection $\sigma =(\cO_{\PP^n},
\cO_{\PP^n}(1),\cdots ,\cO_{\PP^n}(n))$ of locally free sheaves on
$\PP^n$ and we denote by $H$ a hyperplane section on $\PP^n$.  For
any $n \geq r \geq 0$, we have the exact sequence \[0 \rightarrow
\cO_{\PP^n}(r-1) \rightarrow \cO_{\PP^n}(r) \rightarrow \cO_H(r)
\rightarrow 0 .\] According to Proposition \ref{regEi},
$Reg_{\sigma}(\cO_{\PP^n}(r))=-r$ and
$Reg_{\sigma}(\cO_{\PP^n}(r-1))=-r+1$. So,
$Reg_{\sigma}(\cO_{\PP^n}(r)) < max
\{Reg_{\sigma}(\cO_{\PP^n}(r-1)), Reg_{\sigma}(\cO_H(r)) \}$.
\end{example}

\vspace{3mm}

 In Proposition \ref{CM}, we have seen that a coherent sheaf $\cF $ on $\PP^n$ is
 $m$-regular in the sense of Castelnuovo-Mumford if and only if it is
 $m$-regular with respect to the geometric collection
$\sigma =(\cO_{\PP^n}, \cO_{\PP^n}(1),\cdots ,\cO_{\PP^n}(n))$. A
more subtle problem is
 the following one. Let $X \hookright{i} \PP^n$ be a smooth
 projective variety, let $\cF$ be a coherent sheaf on $X$ and let
 $\sigma=(\cE_0, \cdots, \cE_n)$ be a geometric collection of
 coherent sheaves on $X$. How are related the regularity of $\cF$
 with respect to $\sigma$  and
 the Castelnuovo-Mumford regularity of its extension by zero $i_*\cF$ via the embedding
$X \hookright{i} \PP^n$ ?. In next section, we will address this
problem for the case of coherent sheaves on odd dimensional
quadric hypersurfaces $Q_n \subset \PP^{n+1}$.


\section{Regularity of sheaves on $Q_{n}$}

Let $n\in \ZZ$ be an odd integer ($n=2t+1$) and let $Q_n\subset
\PP^{n+1}$ be a smooth quadric hypersurface. In \cite{Kap2}, M. M.
Kapranov defined the locally free sheaves $\psi _i$, $i\ge 0$, on
$Q_n$ and the Spinor bundle $\Sigma $ on $Q_n$ to construct a
resolution of the diagonal $\Delta \subset Q_n \times Q_n$ and to
describe the bounded derived category $D^b({\cO}_{Q_n}-mod)$. In
particular, he got that
\[(\Sigma (-n), \cO_{Q_n}(-n+1), \cdots, \cO_{Q_n}(-1), \cO_{Q_n}) \]
is a full strongly exceptional collection of locally free sheaves
on $Q_n$ (\cite{Kap2}; Proposition 4.9) and hence, according to
Definition \ref{defexcellent}, $Q_n$ has a geometric collection of
locally free sheaves. Dualising each bundle of the above geometric
sequence and reversing the order, we get that
\[\sigma:= (\cO_{Q_n}, \cO_{Q_n}(1), \cdots,  \cO_{Q_n}(n-1), \Sigma (n-1)) \]
is also a geometric collection of locally free sheaves on $Q_n$.

In this section we will give an elementary description of the
locally free sheaves $\psi _i$ and their basic properties (for
more details the reader can look at \cite{Kap2}) and, for any
coherent sheaf $\cF $ on $Q_n$, we will relate $Reg_{\sigma }(\cF
)$ to the Castelnuovo-Mumford regularity of its extension by zero
$i_* \cF$  in the embedding $i:Q_n\hookrightarrow \PP^{n+1}$.

\vskip 2mm From now on, we set $\Omega ^{j}:= \Omega
^j_{\PP^{n+1}}$ and we define inductively $\psi _j$: $$\psi
_0:=\cO_ {Q_n}, \quad \psi _1:=\Omega^1(1)_{|Q_n}$$ and, for all
$j\ge 2$, we define the locally free sheaf $\psi _j$ as the unique
non-splitting extension (Note that $\Ext^1(\psi _{j-2},\Omega
^{j}(j)_{|Q_n})=\CC$): $$0\longrightarrow \Omega
^{j}(j)_{|Q_n}\longrightarrow \psi _j\longrightarrow \psi _{j-2}
\longrightarrow 0. $$

In particular, $\psi _{j+2}=\psi _{j}$ for $j\ge n$ and $\psi
_n=\Sigma (-1)^{2^{t+1}}$, ($n=2t+1$).

\vspace{3mm}

Before computing the right dual basis of any thread
$\sigma_j=(\cE_j, \cE_{j+1}, \cdots, \cE_{j+n})$ of the helix
$\cH_{\sigma}= \{\cE_j \}_{j \in \ZZ}$ associated to $\sigma$, we
collect in the following Lemma the cohomological properties of the
Spinor bundles we need later.

\vspace{3mm}

\begin{lemma} \label{tecnicquadrica} Let $n\in \ZZ$ be an odd integer, let $Q_n\subset
\PP^{n+1}$ be a smooth quadric hypersurface and let $\Sigma$ be
the Spinor bundle on $Q_n$. Then,

(i) $H^i(Q_n,\Sigma(t))=0$ for any $i$ such that $0 < i < n$ and
for all $t \in \ZZ$.

(ii) $H^0(Q_n,\Sigma(t))=0$ for all $t <0$ and $h^0(Q_n,
\Sigma)=2^{\frac{n+1}{2}} $.

(iii) $\Ext^i(\Sigma(j), \Sigma)=\begin{cases} 0 & \mbox{if} \quad
i \neq j \\ \CC & \mbox{if} \quad i= j\end{cases}$ \hspace{3mm}
for any $1 \leq j \leq n$.

(iv) $H^i(Q_n,\psi_j(-j))= \begin{cases} 0 & \mbox{if} \quad i
\neq j \\ \CC & \mbox{if} \quad i= j\end{cases}$ \hspace{3mm} for
any $1 \leq j \leq n$.
\end{lemma}
\begin{proof}
The assertions $(i)$ and $(ii)$ follows from \cite{Ott}; Theorem
2.3. The assertion $(iii)$ follows by induction on $j$, using
$(i)$ and $(ii)$ and the long exact sequence \[ \cdots \rightarrow
\Ext^{i-1}(\cO_{Q_n}(j)^{2^{{\frac{n+1}{2}}}}, \Sigma) \rightarrow
\Ext^{i-1}(\Sigma(j-1), \Sigma) \rightarrow \Ext^i(\Sigma(j),
\Sigma) \rightarrow \Ext^{i}(\cO_{Q_n}(j)^{2^{{\frac{n+1}{2}}}},
\Sigma) \rightarrow \cdots \] obtained by applying the functor
$\Hom( \cdot, \Sigma)$ to the exact sequence $$ 0 \rightarrow
\Sigma(j-1) \rightarrow \cO_{Q_n}(j)^{2^{{\frac{n+1}{2}}}}
\rightarrow \Sigma(j) \rightarrow 0. $$ The assertion $(iv)$
follows from \cite{Kap2}; Proposition 4.11.
\end{proof}

\vspace{3mm}

\begin{proposition}
\label{totselsduals} Let $n\in \ZZ$ be an odd integer, let
$Q_n\subset \PP^{n+1}$ be a smooth quadric hypersurface and let
$\cH_{\sigma}= \{ \cE_i\}_{i  \in \ZZ}$ be the helix associated to
$$\sigma=(\cO_{Q_n}, \cO_{Q_n}(1), \cdots,  \cO_{Q_n}(n-1), \Sigma
(n-1)).$$ Let us denote by $\sigma_k$ the thread $(\cE_k,
\cE_{k+1}, \cdots, \cE_{k+n})$. Then:

(a) The right dual base of $\sigma_0$ is
\[(\Sigma(n-1), \psi_{n-1}(n), \psi_{n-2}(n), \cdots, \psi_1(n), \psi_0(n)).\]

(b) For any $j$, $1 \leq j \leq n$, the right dual base of the
geometric collection
\[\sigma_j=(\cO_{Q_n}(j),
\cdots,  \cO_{Q_n}(n-1), \Sigma (n-1), \cO_{Q_n}(n),
\cO_{Q_n}(n+1), \cdots, \cO_{Q_n}(n+j-1))\] is
\[ \begin{array}{c} (\cO_{Q_n}(n+j-1), \psi_{1}^*(n+j-1), \cdots,
\psi_{j-1}^*(n+j-1), \\ \Sigma(n+j-1), \psi_{n-j-1}(n+j), \cdots,
\psi_0(n+j)).\end{array} \]

(c) For  any $\lambda \in \ZZ$, the right dual base of the
geometric collection
\[ \sigma_{ \lambda(n+1)}=  (\cO_{Q_n}( \lambda n), \cO_{Q_n}(1+ \lambda n),
\cdots,  \cO_{Q_n}(n-1 + \lambda n), \Sigma ( n-1+ \lambda n)) \]
is
\[(\Sigma(n-1 + \lambda n), \psi_{n-1}((\lambda+1)n), \psi_{n-2}((\lambda+1)n),
\cdots, \psi_1((\lambda+1)n), \psi_0((\lambda+1)n)).\]

(d) For any  $j$, $1 \leq j \leq n$ and any $\lambda \in \ZZ$, the
right dual base of the geometric collection
\[\begin{array}{ll} \sigma_{j+ \lambda(n+1)}= & (\cO_{Q_n}(j+ \lambda n),
\cdots,  \cO_{Q_n}(n-1 + \lambda n), \Sigma ( n-1+\lambda n),
\\ &
 \cO_{Q_n}((\lambda+1) n),  \cdots,
\cO_{Q_n}(j-2 + (\lambda+1) n),  \cO_{Q_n}(j-1 + (\lambda+1) n))
\end{array} \] is
\[ \begin{array}{c}(\cO_{Q_n}((\lambda+1)n+j-1), \psi_{1}^*((\lambda+1)n+j-1),
 \cdots,
\psi_{j-1}^*((\lambda+1)n+j-1),   \\   \Sigma((\lambda+1)n+j-1),
\psi_{n-j-1}((\lambda+1)n+j), \cdots, \psi_0((\lambda+1)n+j)).
\end{array}\]
\end{proposition}
\begin{proof}

(a) By \cite{Kap2}; Proposition 4.11 and using the fact that the
right dual basis of an exceptional collection  is uniquely
determined up to unique isomorphism by the orthogonality
conditions described in Remark \ref{dual}, we get that the right
dual basis of the exceptional collection $ (\Sigma (-n),
\cO_{Q_n}(-n+1), \cdots, \cO_{Q_n}(-1), \cO_{Q_n})$ is $$
 (\cO_{Q_n}, \psi_1^*, \cdots, \psi_{n-1}^*, \Sigma^*(1)).
$$
In particular, we have
\begin{equation}
 \label{dualkap} R^{(j)}\cO_{Q_n}(-j)=\psi_j^* \mbox{ for } 0\le j \le n-1; \mbox{ and }
  \end{equation}
$$ R^{(n)}\Sigma (-n)=\Sigma ^*(1).$$

 Let us now compute
the right dual base of $\sigma_0$. Since $R^{(0)}
\Sigma(n-1)=\Sigma(n-1)$ and $R^{(n)} \cO_{Q_n}= \cO_{Q_n}(n)$,
according to the definition of right dual base we only need to
compute $R^{(j)} \cO(n-j)$ for $1 \leq j \leq n-1$. It follows
from Lemma \ref{tecnic}  and the fact that $\Sigma^*(-n+1)=
\Sigma(-n)$, that
\[ \begin{array}{rl} R^{(j)} \cO(n-j)= & R_{\Sigma(n-1) \cO_{Q_n}(n-1) \cdots \cO_{Q_n}(n-j+1)}\cO(n-j) \\ =&
(L_{\Sigma(-n) \cO_{Q_n}(-n+1) \cdots \cO_{Q_n}(j-n-1)}\cO(j-n))^*
\\ = &  (L_{\Sigma \cO_{Q_n}(1) \cdots \cO_{Q_n}(j-1)}\cO(j))^*
\otimes \cO(n).\end{array} \]
By \cite{GK}; 2.8.1,
\[ \begin{array}{rl}(L_{\Sigma
\cO_{Q_n}(1) \cdots \cO_{Q_n}(j-1)}\cO(j))^* \otimes \cO(n)= &
(R_{\cO_{Q_n}(-n+j+1), \cO_{Q_n}(-n+j+2) \cdots
\cO_{Q_n}}\cO_{Q_n}(-n+j))^* \otimes \cO_{Q_n}(n) \\ = &
 (R^{(j)}\cO_{Q_n}(-j))^*
 \otimes \cO_{Q_n}(n)  \\ = &\psi_{n-j}(n) \end{array} \] where the last equality follows
 from (\ref{dualkap}). Hence,
$R^{(j)} \cO(n-j)= \psi_{n-j}(n)$ which finishes the proof of
$(a)$.

(b) Applying Lemma \ref{tecnic}; (c) to (\ref{dualkap}), we get
that
\begin{equation} \label{estrella} R^{(p)} \cO_{Q_n}(n+j-1-p)= \psi_p^*(n+j-1) \quad \mbox{for} \quad
0 \leq p \leq j-1. \end{equation} On the other hand, by Lemma
\ref{tecnicquadrica} we have:

 - For any $\alpha \in \ZZ$ and $ j \leq t \leq n+j-1$
\[ \begin{array}{l}
  \Ext^{\alpha}(\Sigma(n+j-1), \cO_{Q_n}(t))= 0;  \quad \mbox{and} \\
  \Ext^{\alpha}(\Sigma(n+j-1), \Sigma(n-1))= \begin{cases} 0  & \alpha \neq
j\\ \CC & \alpha=j.
 \end{cases}  \\ \end{array} \]

- For any  $\alpha \in \ZZ$, $n-j-1  \leq t  \leq 0$ and $j \leq
\gamma \leq  n+j-1$
\[  \begin{array}{l}
 \Ext^{\alpha}(\psi_{t}(n+j), \Sigma(n-1))=0; \quad \mbox{and}
 \\ \Ext^{\alpha}(\psi_{t}(n+j), \cO_{Q_n}(\gamma))= \begin{cases} \CC
\quad \mbox{if} \quad  n-\alpha=t=\gamma-j  \\ 0 \quad
\mbox{otherwise}.
\end{cases} \\  \end{array}\]

 Since the right dual base of an exceptional
collection is uniquely determined up to unique isomorphism by the
orthogonal conditions given in Remark \ref{dual}, it follows from
(\ref{estrella}), these last cohomological relations and from
Lemma \ref{tecnicquadrica} that indeed
\[ \begin{array}{c}(\cO_{Q_n}(n+j-1), \psi_{1}^*(n+j-1), \psi_{2}^*(n+j-1), \cdots,
\psi_{j-1}^*(n+j-1),  \\ \Sigma(n+j-1), \psi_{n-j-1}(n+j), \cdots,
\psi_1(n+j), \psi_0(n+j)) \end{array} \] is the right dual base of
$\sigma_j$.

$(c)$ and $(d)$ follows from $(a)$ and $(b)$, respectively,  and
Corollary \ref{corotecnic}.
\end{proof}

Let $\cF$ be a coherent sheaf on $Q_n$ and let $i_*\cF$ be its
extension by zero in the embedding $Q_n \hookright{i} \PP^{n+1}$.
We are now ready to compare $Reg_{\sigma}(\cF)$ to
$Reg^{CM}(i_*\cF)$ and state the main result of this section.

\begin{theorem}
\label{mainquadrica} Let $n$ be an odd integer and let $Q_n
\hookright{i} \PP^{n+1}$ be a quadric hypersurface. We consider
the geometric collection  $\sigma=(\cO_{Q_n}, \cdots,
\cO_{Q_n}(n-1), \Sigma (n-1))$. For any coherent sheaf $\cF$ on
$Q_n$ we have:
  \[ \lfloor \frac{nReg_{\sigma}(\cF)}{n+1}  \rfloor  \leq
  Reg^{CM}(i_*\cF) \leq  \lfloor \frac{nReg_{\sigma}(\cF)}{n+1} \rfloor  +1 .\]

\end{theorem}
\begin{proof}
First of all, we will see that $Reg^{CM}(i_*\cF) \leq \lfloor
\frac{nReg_{\sigma}(\cF)}{n+1}  \rfloor  +1$. To this end, since
for all $t \in \ZZ$ and $0 \leq q \leq n$, \begin{equation}
\label{pnQn} H^{n+1}(\PP^{n+1}, i_*\cF(t))=0, \quad \quad
H^q(\PP^{n+1}, i_*\cF(t)) \cong H^q(Q_{n}, \cF(t)),
\end{equation}
we will see that if $\cF$ is $m$-regular with respect to $\sigma$
then
\[ H^i(Q_n, \cF(m+\lambda+1-i))=0, \quad 1 \leq i \leq n \]
where $\lambda:= \lfloor \frac{-m}{n+1}\rfloor$. Write
$-m=\lambda(n+1)+r$ with $0 \leq r \leq n$. So,
\[-m-n= \begin{cases} (\lambda-1)(n+1)+r+1 &  \mbox{if } 0 \leq r \leq n-1 \\
  \lambda(n+1) & \mbox{if }  r=n. \\
\end{cases} \]
Therefore, $\sigma_{-m-n}=( \cE_{-m-n}, \cdots, \cE_{-m})$ is
equal to
\[ \begin{array}{c}(\cO_{Q_n}((\lambda-1)n+r+1), \cdots,  \cO_{Q_n}((\lambda-1)n+n-1),
\Sigma((\lambda-1)n+n-1),   \\  \cO_{Q_n}(\lambda n ),
\cO_{Q_n}(\lambda n +1), \cdots,  \cO_{Q_n}({\lambda n +r})),
\end{array}\] if $0 \leq r \leq n-1$ and is equal to \[ (\cO_{Q_n}(\lambda n),
\cO_{Q_n}(\lambda n +1), \dots, \cO_{Q_n}(\lambda n + n-1),
\Sigma(\lambda n +n-1) ) \] if $r=n$.  We consider the
Eilenberg-Moore type spectral sequence (see Theorem \ref{moore})

$$ E_1^{pq}=\bigoplus _{\alpha + \beta =q} \Big(\Ext^{\alpha
}(R^{(-p)}\cE _{i+n+p},\cF )\otimes \Ext^{\beta }(\cG,\cE
_{i+n+p})\Big) \quad 0\le q \le n, \quad -n\le p \le 0$$

$$d_r^{pq}:E_r^{pq} \longrightarrow E_{r}^{p+r,q-r+1}, \quad \quad
E_{\infty }\Rightarrow \Ext^{p+q}(\cG,\cF) $$ and we apply it to
the case $i=-m-n$ and $\cG=\cO _{Q_n}(t-m- \lambda-1)$. So, in the
square $ 0\le q \le n$, $-n\le p \le 0$ we have $$
E_1^{pq}=\bigoplus _{\alpha + \beta =q} \Big(\Ext^{\alpha
}(R^{(-p)}\cE _{-m+p},\cF )\otimes \Ext^{\beta }(\cO _{Q_n}(t-m-
\lambda-1),\cE _{-m+p})\Big). $$ Since $\cF$ is $m$-regular with
respect to $\sigma$, $\Ext^{\alpha }(R^{(-p)}\cE _{-m+p},\cF )=0$
for all $\alpha>0$ and all $p$, $-n \leq p \leq 0$. On the other
hand, if $0 \leq r \leq n-1$, then \[ \begin{array}{ll}
\Ext^{\beta }(\cO _{Q_n}(t-m- \lambda-1),\cE _{-m+p})  & =
\begin{cases} H^{\beta}(Q_n, \cO_{Q_n}(2-t+p)) & \mbox{if } -n
\leq p \leq -r-2
\\ H^{\beta}(Q_n, \Sigma(-t+p+1)) & \mbox{if }  p = -r-1
\\
H^{\beta}(Q_n, \cO_{Q_n}(1-t+p)) & \mbox{if } -r \leq p \leq 0,
\end{cases} \end{array}\]
and if $r=n$, then
 \[ \begin{array}{ll} \Ext^{\beta }(\cO _{Q_n}(t-m-
\lambda-1),\cE _{-m+p}) & =
\begin{cases} H^{\beta}(Q_n, \cO_{Q_n}(1-t+p)) & \mbox{if } -n
\leq p \leq -1
\\ H^{\beta}(Q_n, \Sigma(-t)) & \mbox{if }  p = 0.
\end{cases} \end{array}\]
Therefore, applying Serre's duality and Lemma
\ref{tecnicquadrica}, we get that the only non-zero Ext's groups
$\Ext^{\beta}(\cO _{Q_n}(t-m- \lambda-1),\cE _{-m+p}))$ correspond
to $\beta=0$ or $\beta=n$ and in the latter case $-r \leq p \leq
0$ and $p \leq t-1-n$ or $-n \leq p \leq -r-1$ and $p \leq t-2-n$.
Thus, $E_1^{pq} \neq 0$ forces $p+q \leq t-1$. Therefore,
$E_{\infty}^{pq}=0$ for $p+q \geq t$ and so for any $t
>0$, $$\Ext^{t}(\cO _{Q_n}(t-m- \lambda-1), \cF)= H^t(Q_n, \cF(m +
\lambda+1-t))=0.$$

In order to prove the other inequality, we will see that if
$i_*\cF$ is $(m+ \lambda)$-regular in the sense of
Castelnuovo-Mumford, then $\cF$ is $m$-regular with respect to
$\sigma$, where $\lambda:= \lfloor \frac{-m}{n+1}\rfloor$. If
$i_*\cF$ is $m$-regular, then it is also $(m+\gamma)$-regular for
all $\gamma \geq 0$ and we have
\begin{equation}
\label{condicio} 0= H^i(\PP^{n+1}, i_*\cF(t))=H^i(Q_n, \cF(t)),
\quad \mbox{for all } t \geq m+\lambda-i, \quad i>0.
\end{equation}

We write $-m=\lambda (n+1)+r$ with $0 \leq r \leq n$. If $r \neq
n$ then $\sigma_{-m-n}$ is given by
\[ \begin{array}{c}(\cO_{Q_n}((\lambda-1)n+r+1), \cdots,  \cO_{Q_n}((\lambda-1)n+n-1),
\Sigma((\lambda-1)n+n-1),   \\  \cO_{Q_n}(\lambda n ),
\cO_{Q_n}(\lambda n +1), \cdots,  \cO_{Q_n}({\lambda n +r})),
\end{array}\]
 and by Proposition \ref{totselsduals} its right dual base is given by
\[ \begin{array}{c}(\cO_{Q_n}(\lambda n+r), \psi_{1}^*(\lambda n+r),
 \cdots,
\psi_{r}^*(\lambda n+r),   \\   \Sigma(\lambda n+r),
\psi_{n-r-2}(\lambda n+r+1), \cdots, \psi_0(\lambda n+r+1)),
\end{array}\]
and if $r=n$,   $\sigma_{-m-n}$ is equal to \[ (\cO_{Q_n}(\lambda
n), \cO_{Q_n}(\lambda n +1), \dots, \cO_{Q_n}(\lambda n + n-1),
\Sigma(\lambda n +n-1) ) \] and by Proposition \ref{totselsduals}
its right dual base is given by
\[(\Sigma(n-1 + \lambda n), \psi_{n-1}((\lambda+1)n), \psi_{n-2}((\lambda+1)n),
\cdots, \psi_1((\lambda+1)n), \psi_0((\lambda+1)n)).\] So,
according to Definition \ref{new}, we have to see that for any $q
>0$
\[  \begin{array}{lll}
(i) & \Ext^q(\cO_{Q_n}(\lambda n +r), \cF)=H^q(Q_n, \cF(m+
\lambda))=0 & \\
 (ii) & \Ext^q(\psi_{\alpha}^*(\lambda n +r),\cF)=
 \Ext^q(\psi_{\alpha}^*, \cF(m+ \lambda))=0 & 1 \leq \alpha \leq r \\
 (iii) & \Ext^q(\psi_{n-\beta}(\lambda n +r+1),\cF)=
 \Ext^q(\psi_{n-\beta}, \cF(m+ \lambda-1))=0 & r+2 \leq \beta \leq n \\
 (iv) & \Ext^q(\Sigma(\lambda n +r),\cF)=
 \Ext^q(\Sigma, \cF(m+ \lambda))=0 &  \\
\end{array}\]
if $r \neq n$ and we have to see that
\[  \begin{array}{lll}
 (v) & \Ext^q(\psi_{n-\alpha}((\lambda+1) n ),\cF)=
 \Ext^q(\psi_{n-\alpha}, \cF(m+ \lambda))=0 & 1 \leq \alpha \leq n \\
 (vi) & \Ext^q(\Sigma(\lambda n +n-1),\cF)=
 \Ext^q(\Sigma, \cF(m+ \lambda+1))=0 &  \\
\end{array}\] if $r=n$. From (\ref{condicio}) we immediately get $(i)$.
Equalities $(ii)$, $(iii)$ and $(v)$ are consequence of the
following stronger results
\begin{equation}
\label{resultat} \Ext^{\alpha}(\psi_j, \cF(t))=0  \quad \mbox{for
any } \alpha>0, \quad 0 \leq j \leq n, \quad \mbox{and } t \geq
m+\lambda-\alpha,
\end{equation}
and
\begin{equation}
\label{resultat2} \Ext^{\alpha}(\psi_j^*, \cF(t))=0  \quad
\mbox{for any } \alpha>0, \quad 1 \leq j \leq r, \quad \mbox{and }
t \geq m+\lambda-\alpha.
\end{equation}

We will prove (\ref{resultat}) by induction on $j$. Analogously,
we can prove (\ref{resultat2}) and we right it to the reader. If
$j=0$, it follows from (\ref{condicio}) that
$\Ext^{\alpha}(\psi_0, \cF(t))=H^{\alpha}(Q_n, \cF(t))=0$ for all
$\alpha>0$ and $t \geq m+\lambda-\alpha$. Assume $j=1$.
Restricting to $Q_n$ the Euler sequence
 and applying  the functor $\Ext(\cdot, \cF(t))$ we get the
 exact sequence
 \[ \cdots \rightarrow \Ext^{\alpha}(\cO_{Q_n}^{n+2}, \cF(t)) \rightarrow \Ext^{\alpha}(\psi_1, \cF(t))
 \rightarrow \Ext^{\alpha+1}(\cO_{Q_n}(1), \cF(t)) \rightarrow \cdots . \]
 By (\ref{condicio}),   $\Ext^{\alpha+1}(\cO_{Q_n}(1), \cF(t))=
 H^{\alpha+1}(Q_n, \cF(t-1))=0 $ for $t -1 \geq m + \lambda
 -\alpha-1$ and $\Ext^{\alpha}(\cO_{Q_n}^{n+2}, \cF(t))=0$ for all
 $\alpha >0$ and $t \geq m + \lambda - \alpha$. Hence
 \begin{equation}
 \label{j1}
 \Ext^{\alpha}(\psi_1, \cF(t))=0 \quad \mbox{for all } \alpha>0
 \quad \mbox{and } t \geq m + \lambda - \alpha.
 \end{equation}

 Before going ahead with the general case, let us prove the
 following
\vskip 2mm

\noindent
 {\bf Claim:} For any $p \geq 1$, $\alpha>0$ and $t \geq m + \lambda -
 \alpha$,
 \[ \Ext^{\alpha}(\Omega^p(p)_{|_{Q_n}}, \cF(t))=0.\]
 {\bf Proof of the Claim:} We will prove it by induction on $p$.
 If $p=1$, by (\ref{j1}),
 \[ \Ext^{\alpha}(\Omega^1(1)_{|_{Q_n}}, \cF(t))= \Ext^{\alpha}(\psi_1, \cF(t))=0\]
 for any $\alpha>0$ and $t \geq m + \lambda -
 \alpha$. Assume it holds for $p$ and let us see the case $p+1$.
 Applying the functor $\Ext(\cdot, \cF(t))$ to the exact sequence
\[0 \rightarrow \Omega^{p+1}(p+1)_{|_{Q_n}} \rightarrow \cO_{Q_n}^{n+2 \choose p+1}
 \rightarrow \Omega^{p}(p+1)_{|_{Q_n}} \rightarrow 0 \]
 obtained by restricting to $Q_n$ the $p$-th exterior power of the
 dual of the Euler sequence on $\PP^{n+1}$, we get the long exact
 sequence
  \[ \cdots \rightarrow \Ext^{\alpha}(\cO_{Q_n}^{n+2 \choose p+1}, \cF(t))
   \rightarrow \Ext^{\alpha}(\Omega^{p+1}(p+1)_{|_{Q_n}}, \cF(t))
 \rightarrow \Ext^{\alpha+1}(\Omega^{p}(p+1)_{|_{Q_n}}, \cF(t)) \rightarrow \cdots . \]
 By hypothesis of induction, $\Ext^{\alpha+1}(\Omega^{p}(p+1)_{|_{Q_n}},
 \cF(t))=0$ for any $t-1 \geq m + \lambda - \alpha -1$ and by
 (\ref{condicio}), $\Ext^{\alpha}(\cO_{Q_n}^{n+2 \choose p+1},
 \cF(t))=0$ for any $t \geq m + \lambda - \alpha$. Hence,
  $\Ext^{\alpha}(\Omega^{p+1}(p+1)_{|_{Q_n}}, \cF(t))=0$ for any
  $t \geq m + \lambda - \alpha$ and $\alpha>0$ which finishes the
  proof of the Claim.

\vspace{3mm}

  Let us now prove (\ref{resultat}) by induction on $j$. The first
  two cases have been already done. Assume (\ref{resultat}) holds
  for $j \geq 1$ and let us see that it holds for $j+1$.
  Applying the functor $\Ext(\cdot, \cF(t))$ to the exact sequence
$$0\longrightarrow \Omega ^{j}(j)_{|Q_n}\longrightarrow \psi
_j\longrightarrow \psi _{j-2} \longrightarrow 0 $$ we get the long
exact sequence
\[ \cdots \rightarrow \Ext^{\alpha}(\psi_{j-1}, \cF(t)) \rightarrow
\Ext^{\alpha}(\psi_{j+1}, \cF(t)) \rightarrow
\Ext^{\alpha}(\Omega^{j+1}(j+1)_{|_{Q_n}}, \cF(t)) \rightarrow
\cdots .\] It follows from the Claim that
$\Ext^{\alpha}(\Omega^{j+1}(j+1)|_{Q_n}, \cF(t))=0$ for all $
\alpha>0$ and $t \geq m + \lambda - \alpha$ and by hypothesis of
induction $\Ext^{\alpha}(\psi_{j-1}, \cF(t))=0$ for all $
\alpha>0$ and $t \geq m + \lambda - \alpha$. Hence we obtain
\[ \Ext^{\alpha}(\psi_{j+1}, \cF(t))=0\]
for any $\alpha >0$ and $t \geq m + \lambda - \alpha$, which
finishes the proof of $(ii)$, $(iii)$ and $(v)$. Finally, using
the fact that $\psi_n= \Sigma(-1)^{2^{t+1}}$, we deduce from
(\ref{resultat}) that $(iv)$ and $(vi)$ also hold.
\end{proof}

Now we will show that Theorem \ref{mainquadrica} is optimal, i.e.,
both inequalities can be realized. First,  we will compare
$Reg_{\sigma}(\cE_j)$ to $Reg^{CM}(i_*\cE_j)$ where $\cE_j$ is
part of the helix $\cH_{\sigma}= \{ \cE_j\}_{j \in \ZZ}$
associated to $\sigma$ and later we will compare
$Reg_{\sigma}(\psi_1(3+\lambda n))$ to $Reg^{CM}(i_*(\psi_1(3 +
\lambda n )))$.

\begin{proposition}
\label{regEiquadrica} Let $n$ be an odd integer and let $i: Q_n
\hookrightarrow \PP^{n+1}$ be a quadric hypersurface. Let
$\cH_{\sigma}= \{ \cE_j\}_{j \in \ZZ}$ be the strict helix
associated to $\sigma=(\cO_{Q_n}, \cdots, \cO_{Q_n}(n-1), \Sigma
(n-1))$. For any integer $j\in \ZZ$, we write $j=\lambda(n+1)+r$
with  $0 \leq r \leq n$. Then, we have
\[ \begin{array}{l} Reg_{\sigma}(\cE_j)=-j \\
  Reg^{CM}(i_* \cE_j)=-j+\lambda+1. \end{array}\]
  In particular, the following relation holds:
  \[Reg_{\sigma}(\cE_j)+ \lambda+1= Reg^{CM}(i_*\cE_j). \]
\end{proposition}
 \begin{proof}
 By remark \ref{period}, if $0 \leq r \leq n-1$, we have
 \[ \cE_j=\cE_{\lambda(n+1)+r}= \cE_r \otimes (K^*)^{\otimes \lambda} = \cO_{Q_n}(\lambda n +r)\]
 and if $r=n$, then we have
 \[ \cE_j=\cE_{\lambda(n+1)+n}= \cE_n \otimes (K^*)^{\otimes \lambda} = \Sigma(\lambda n +n-1).\]
 Using the exact sequence
 \[ 0 \rightarrow \cO_{Q_n}(-2) \rightarrow \cO_{Q_n} \rightarrow i_*\cO_{Q_n} \rightarrow 0\]
 we get that
 \[ \begin{array}{ll} H^q(\PP^{n+1}, i_*\cO_{Q_n}(t))=  0 & \mbox{for all } t \in \ZZ,
 \quad 1 \leq q \leq n-1 \quad
  \mbox{and } q=n+1 \\  H^n(\PP^{n+1}, i_*\cO_{Q_n}(t))=     0 & \mbox{for all
  } t \geq -n+1,
    \end{array} \]
    and we conclude that $Reg^{CM}(i_*\cO_{Q_n}(s))=-s+1$. In
    particular, if $j= \lambda(n+1)+r$ with $0 \leq  r \leq n-1$,
    we have
    \[ \begin{array}{ll}Reg^{CM}(i_*\cE_j) & = Reg^{CM}( i_*\cO_{Q_n}(\lambda n +r))
    \\ & = -\lambda n -r +1 \\ &= -j + \lambda +1 \\
    & =     Reg_{\sigma}(\cE_j) + \lambda +1 \end{array}\]
    where the last equality follows from Proposition \ref{regEi}.
    So, it only remains the case $j= \lambda(n+1)+n$. Using the
    isomorphism
    \[ H^q(\PP^{n+1}, i_* \Sigma(t)) \cong H^q(Q_n, \Sigma(t))\]
    and Lemma \ref{tecnicquadrica}; $(i)-(ii)$, we obtain that $Reg^{CM}(i_*
    \Sigma(t))=-t$. Therefore, if $j= \lambda(n+1)+n$,
    we have
       \[\begin{array}{ll}Reg^{CM}(i_*\cE_j) & = Reg^{CM}( i_*\Sigma(\lambda n +n-1)) \\ &
        = -\lambda n -n +1 \\ & =
        -j + \lambda +1 \\ & =
    Reg_{\sigma}(\cE_j) + \lambda +1 \end{array} \]
    where again the last equality follows from Proposition
    \ref{regEi}.
 \end{proof}

\begin{proposition}
\label{exemplepsi1} With the above notations, we have:
\[ \begin{array}{ll} (i) & Reg^{CM}(i_*(\psi_1(3 + \lambda
n)))=-2-\lambda n. \\ (ii) & Reg_{\sigma}(\psi_1(3 + \lambda
n))=-2-\lambda (n+1).
\end{array} \]
\end{proposition}
\begin{proof} $(i)$ It follows from the fact that $H ^
i(Q_n,\psi_1(1-i))=0$ for all $i \geq 1$ and
\[ H^1(Q_n,\psi_1(-1)) \cong \CC. \]

$(ii)$ Since by Proposition \ref{totselsduals} the right dual
basis of $\sigma_{2+\lambda(n+1)-n}$ is
\[ (\cO_{Q_n}(\lambda n +2), \psi_{1}^*(\lambda n +2), \psi_{2}^*(\lambda n +2),
\Sigma(\lambda n +2), \psi_{n-4}(\lambda n +3), \cdots,
\psi_{0}(\lambda n +3))
\]
and the right dual basis of $\sigma_{3+\lambda(n+1)-n}$ is
\[ (\cO_{Q_n}(\lambda n +3), \psi_{1}^*(\lambda n +3), \psi_{2}^*(\lambda n +3), \psi_{3}^*(\lambda n +3),
\Sigma(\lambda n +3), \psi_{n-5}(\lambda n +4), \cdots,
\psi_{0}(\lambda n +4))
\]
we easily check (taking into account that $\psi_0(\lambda n +4)
\cong \cO_{Q_n}(\lambda n +4)$ and that $H^1(Q_n, \psi_1(-1))=
\CC$) that $\psi_1(3 + \lambda n)$ is $(-2-\lambda (n+1))$-regular
with respect to $\sigma$ but not $(-3-\lambda (n+1))$-regular with
respect to $\sigma$.
\end{proof}


\section{Final comments}

\rm In this section, we will first gather the problems that arose
from this paper and we will end with some final remarks.

\begin{problem}\label{prob1}
To characterize smooth projective varieties which have a geometric
collection.
\end{problem}

Let $X$ be a smooth projective variety of dimension $n$. We know
that a necessary condition for having a geometric collection is
that the rank of the Grothendieck group $K_0(X)$ is $n+1$ and $X$
is Fano. We would like to know if these conditions are also
sufficient.

 \rm The main goal of this paper was to generalize the
notion of Castelnuovo-Mumford regularity for coherent sheaves on a
projective space to  coherent sheaves on other smooth projective
varieties $X$. We have succeed provided $X$ has a geometric
collection of coherent sheaves. It would be nice to extend the
definition to the case  of smooth projective varieties which have
a full strongly exceptional collection of sheaves (not necessarily
geometric). Hence, we propose:

\vspace{3mm}

\begin{problem} Let $X$ be a smooth projective variety and let
$\cF$ be a coherent sheaf on $X$. To extend the definition of
regularity of $\cF$ with respect to a geometric collection
$\sigma$ to regularity of $\cF$ with respect to a full strongly
exceptional collection.
\end{problem}

\vskip 2mm  It is well known that Beilinson's theorem (\cite{Be})
and Castelnuovo-Mumford regularity of sheaves play a fundamental
role in the classification of vector bundles on projective spaces.
In a forthcoming paper we will apply the results obtained in this
work to study moduli spaces of vector bundles on quadric
hypersurfaces. We have the feeling that Beilinson-Kapranov type
spectral sequence (Theorem \ref{beil}), Eilemberg-Moore type
spectral sequence (Theorem \ref{moore}) and the regularity with
respect to a geometric collection will play an important role in
the classification of vector bundles on varieties with a geometric
collection.

\end{document}